\documentclass{article}
\usepackage{amsmath,amssymb,amsfonts}
\usepackage{graphicx}
\usepackage{color}

\def\rank{\mathrm{rank}}

\def\Jac{\mathrm{Jac}}
\def\tr{\mathrm{tr}}

\newtheorem{theorem}{Theorem}
\newtheorem{proposition}{Proposition}
\newtheorem{corollary}{Corollary}
\newtheorem{lemma}{Lemma}
\newtheorem{remark}{Remark}
\newtheorem{example}{Example}
\newtheorem{definition}{Definition}
\newenvironment{proof}[1][Proof]{\noindent\textit{#1.} }{\hfill$\Box$\medskip}

\title{Hyperelliptic Jacobians as Billiard Algebra of Pencils of Quadrics: Beyond Poncelet Porisms}

\author{Vladimir Dragovi\'c and Milena Radnovi\'c}

\date{}

\begin{document}

\maketitle


\medskip

\centerline{Mathematical Institute SANU}

\smallskip

\centerline{Kneza Mihaila 36, 11000 Belgrade, Serbia}

\smallskip

\centerline{e-mail: {\tt vladad@mi.sanu.ac.yu,
milena@mi.sanu.ac.yu}}

\vfill

\centerline{The corresponding author:}

\smallskip

\centerline{Vladimir Dragovi\'c}

\smallskip

\centerline{Mathematical Institute SANU}

\centerline{Kneza Mihaila 36, 11000 Belgrade, Serbia}

\smallskip

\centerline{e-mail: {\tt vladad@mi.sanu.ac.yu}}

\smallskip

\centerline{phone: +381 11 2630170}

\centerline{fax: +381 11 2186105}

\

\

\begin{abstract}

\smallskip

The thirty years old programme of Griffiths and Harris of
understanding higher-dimensional analogues of Poncelet-type problems
and synthetic approach to higher genera addition theorems has been
settled and completed in this paper. Starting with the observation
of the billiard nature of some classical constructions and
configurations, we construct {\em the billiard algebra}, that is a
group structure on the set $T$ of lines simultaneously tangent to
$d-1$ quadrics from a given confocal family in the $d$-dimensional
Euclidean space. Using this tool, the related results of Reid,
Donagi and Kn\"orrer are further developed, realized and simplified.
We derive a fundamental property of $T$: any two lines from this set
can be obtained from each other by at most $d-1$ billiard
reflections at some quadrics from the confocal family. We introduce
two hierarchies of notions: {\em $s$-skew lines} in $T$ and {\em
$s$-weak Poncelet trajectories}, $s=-1,0,...,d-2$. The
interrelations between billiard dynamics, linear subspaces of
intersections of quadrics and hyperelliptic Jacobians developed in
this paper enabled us to obtain higher-dimensional and higher-genera
generalizations of several classical genus $1$ results: Cayley's
theorem, Weyr's theorem, Griffiths-Harris theorem and Darboux
theorem.
\end{abstract}

\newpage

\tableofcontents

\newpage

\section{Introduction}

\begin{figure}[h]
\centering
\begin{minipage}[b]{0.6\textwidth}
\centering
\includegraphics[width=5cm,height=4cm]{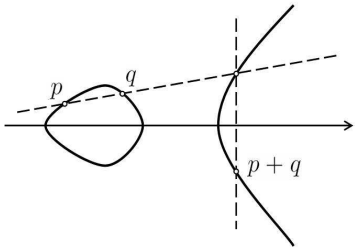}
\caption{The group law on the cubic curve}
\end{minipage}\label{fig:kubika}
\end{figure}

Having in mind the geometric interpretation of the group structure
on a cubic (see Figure \ref{fig:kubika}),
the question of finding an analogous construction of the group
structure in higher genera arises. In this paper, we show that such
a geometric construction exists in the case of hyperelliptic
Jacobians. Our ideas are continuation of those of Reid, Donagi and
Kn\"orrer, see \cite{Re}, \cite{Kn}, \cite{Don} and also \cite{DR},
\cite{N}. Further development, realization, simplification and
visualization of their constructions is obtained by using the ideas
of billiard dynamics on pencils of quadrics developed in \cite{DR4}.

\smallskip

The key projective-geometric feature of that billiard dynamics is
the Double Reflection Theorem, see Theorem \ref{th:DRT} below. There
are four lines belonging to a certain linear space and forming {\it
the Double reflection configuration}: these four lines reflect to
each other according to {\bf the billiard law at some confocal
quadrics}. This billiard configuration appears to be the exact genus
two generalization of three points belonging to a line in the
elliptic case.

\smallskip

In higher genera, we construct the corresponding more general
billiard configuration, again by using Double Reflection Theorem.
This configuration, which we call {\it $s$-brush}, is in one of the
equivalent formulations, a certain billiard trajectory of length
$s\le g$ and the sum of $s$ elements in the brush is, roughly
speaking, the final segment of that billiard trajectory.

\smallskip

The milestones of this paper are \cite{Kn} and \cite{DR4} and the
key observation giving a link between them is {\it that the
correspondence $g\mapsto g'$ in Lemma 4.1 and Corollary 4.2 from
\cite{Kn} is } {\bf the billiard map at the quadric} $\mathcal
Q_{\lambda}$.

\smallskip

Thus, after observing and understanding the billiard nature behind
the constructions of \cite {Re}, \cite{Kn}, \cite{Don}, we become
able to use the billiard tools to construct and study hyperelliptic
Jacobians, and particularly their real part. Any real hyperelliptic
Jacobian may be realized as a set $T$ of lines simultaneously
tangent to given $d-1$ quadrics $\mathcal Q_1$, ..., $\mathcal
Q_{d-1}$ of some confocal family in the $d$-dimensional Euclidean
space. It is well known that such a set $T$ is invariant under the
billiard dynamics determined by quadrics from the confocal family.
By using Double Reflection Theorem and some other billiard
constructions we construct a group structure on $T$, {\it a billiard
algebra}. The usage of billiard dynamics in algebro-geometric
considerations appears to be, as usually in such a situation, of a
two-way benefit. We derive a fundamental property of $T$: {\it any
two lines in $T$ can be obtained from each other by at most $d-1$
billiard reflections at some quadrics from the confocal family}. The
last fact opens a possibility to introduce new hierarchies of
notions: of {\it $s$-skew lines in $T$, $s=-1, 0,\dots, d-2$} and of
{\it $s$-weak Poncelet trajectories of length $n$}. The last are
natural quasi-periodic generalizations of Poncelet polygons. By
using billiard algebra, we obtain complete analytical description
for them. These results are further generalizations of our recent
description of Cayley's type of Poncelet polygons in arbitrary
dimension, see \cite {DR4}. Let us emphasize that the method used
here, based on billiard algebra differs from the methods exposed in
\cite {DR4}.

\smallskip

The interrelations between billiard dynamics, subspaces of
intersections of quadrics and hyperelliptic Jacobians developed in
this paper, enable us to obtain higher dimensional generalizations
of several classical results. To demonstrate the power of our
method, we present here generalizations of Weyr's Poncelet theorem
(see \cite {We}) and also Griffiths-Harris Space Poncelet theorem
(see \cite {GH1}) in arbitrary dimension are derived and presented
here. We also give an arbitrary-dimensional generalization of the
Darboux theorem \cite{Dar3}. Let us mention that one of the plane
versions of Darboux theorem has been recently rediscovered, with
some improvements, in \cite{Svarc}. In that paper, the obtained
configuration is called {\em Poncelet grid}, while we name it here
{\em Poncelet-Darboux grid}, which is, by our opinion, more
historically justified.

\smallskip

The paper is organized as follows. Section \ref{sec:preliminaries}
consists of preliminaries. It starts with a review of confocal
families of quadrics and their most important properties, in the
$d$-dimensional Euclidean and projective space. The Poncelet theorem
in the three-dimensional projective space over an arbitrary field is
proved, giving accent on projective definitions and methods. We
review The One Reflection Theorem and The Double Reflection Theorem
(DRT). The latter appears to be one of the main projective tools
used in the paper. The notions of virtual reflections and virtual
reflection configuration are recalled and some of their basic
properties from \cite{DR4} are reviewed. We recall the definition of
the generalized Cayley's curve from \cite{DR4}, and discuss its form
in the real case.

\smallskip

Section \ref{sekcija:bilijarskaalgebra} is the kernel of the paper.
In Section \ref{sekcija:morfizam}, relying on the simple Lemma
\ref{lema:veza}, we begin to build a bridge between Kn\"orrer's
constructions from \cite{Kn} and the billiard notions from our
previous article \cite{DR4} and this paper. Then, we construct a
birational morphism between the generalized Cayley's curve and the
Reid-Donagi-Kn\"orrer curve (see \cite{Re},\cite{Don},\cite{Kn}). In
Section \ref{sec:genus2}, we construct a group structure on the set
of lines in $\mathbf E^3$ that are simultaneously tangent to two
given confocal quadrics. This group structure is naturally connected
with the billiard law. Using its properties, some beautiful
properties of confocal families are derived, e.g.\ Theorem
\ref{th:zvezda}. In Section \ref{sec:general}, set $T$ of lines
simultaneously tangent to given $d-1$ quadrics $\mathcal Q_1$, ...,
$\mathcal Q_{d-1}$ of some confocal family in $\mathbf E^d$ is in
several steps endowed with a group structure called {\it billiard
algebra}. In this construction, we first define the addition
operation on the set of all finite billiard trajectories with the
same initial line. Then, observing equivalence relations $\alpha$
and $\beta$ in the set of finite billiard trajectories, we prove
that the addition operation is compatible with them, and that the
operation on the equivalence classes has additional algebraic
properties. At the end, we are ready to use these constructions in
order to obtain the billiard algebra on the set of lines.

\smallskip

Section \ref{sec:algebra} is aimed to illustrate the billiard
algebra introduced in Section \ref{sekcija:bilijarskaalgebra}.
Section \ref{sec:weak} starts with introduction of the notions of
{\em $s$-skew lines} and {\em $s$-weak Poncelet trajectories}. We
prove the Poncelet type theorem for such trajectories and deduce the
corresponding analytic condition of Cayley's type. In Section
\ref{sekcija:wgh}, higher-dimensional generalizations of some
classical results connected to Poncelet's porism are stated and
proved. We introduce {\em generalized Weyr chains} and show their
correspondence with Poncelet polygons (Theorem \ref{th:gen.weyr} and
Proposition \ref{propozicija:weyr.poncelet}). This correspondence is
then used to demonstrate the generalization to higher dimensions of
the Griffiths-Harris space version of Poncelet theorem (Theorem
\ref{teorema:uopstena.grifits.haris}). Section \ref{sec:grid} starts
with recalling and a few comments on Darboux's theorem on grids of
conics associated to a Poncelet polygon on an arbitrary Liouville
surface. Then, for the case when the surface is the Eucledean plane,
we prove in Theorems \ref{th:gen.grida} and \ref{th:grida} a
generalization of the Darboux's theorem. The essential
generalization to higher dimensions is demonstrated in Theorem
\ref{th:gen.grida.d}.

\smallskip

Section \ref{sec:conclusion} contains concluding remarks.

\section{Preliminaries}\label{sec:preliminaries}

\subsection{Confocal Families of Quadrics and Billiards in Euclidean
Space}\label{sec:confocal.E}

In this section, we are going to define families of confocal
quadrics in the $d$-dimensional Euclidean space $\mathbf E^d$ and
summarize their basic properties connected with the billiard
reflection law.

\smallskip

\begin{definition}\label{def:E.confocal}
{\em A family of confocal quadrics} in the $d$-dimensional Euclidean
space $\mathbf E^d$ is a family of the form:
\begin{equation}\label{eq:confocal.family}
\mathcal Q_{\lambda}\ :\
\frac{x_1^2}{a_1-\lambda}+\dots+\frac{x_d^2}{a_d-\lambda}=1\qquad(\lambda\in\mathbf
R),
\end{equation}
where $a_1$, \dots, $a_d$ are real constants.
\end{definition}

Let us notice that a family of confocal quadrics in the Euclidean
space is determined by only one quadric.

From now on, we are going to consider the non-degenerate case when
$a_1$,\dots,$a_d$ are all distinct.

\begin{theorem}[Jacobi]\label{th:Jacobi}
Any point of the $d$-dimensional Euclidean space is the intersection
of exactly $d$ quadrics of the confocal family
(\ref{eq:confocal.family}). The quadrics are perpendicular to each
other at the intersecting points. (See Figure \ref{fig:3konf}).
\begin{figure}[h]
\centering
\begin{minipage}[b]{0.6\textwidth}
\centering
\includegraphics[width=6cm,height=6cm]{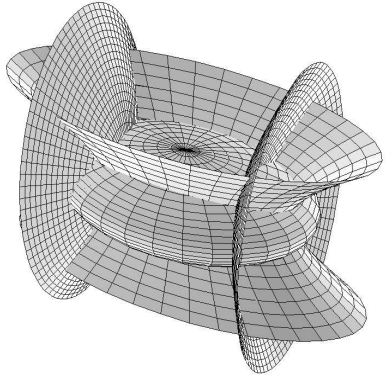}
\caption{Three confocal quadrics in $\mathbf E^3$}\label{fig:3konf}
\end{minipage}
\end{figure}
\end{theorem}

From here, it follows that to every point $x\in\mathbf E^d$ we may
associate a $d$-tuple of distinct parameters
$(\lambda_1,\dots,\lambda_d)$, such that $x$ belongs to quadrics
$\mathcal Q_{\lambda_1}$,\dots, $\mathcal Q_{\lambda_d}$.
Additionally, if we order $a_1<a_2<\dots<a_d$ and
$\lambda_1<\lambda_2\dots<\lambda_d$, then:
$$
\lambda_1\le a_1<\lambda_2\le a_2<\dots<\lambda_d\le a_d.
$$

\begin{definition}
{\em Jacobi elliptic coordinates} of point
$x=(x_1,\dots,x_d)\in\mathbf E^d$ are values
$(\lambda_1,\dots,\lambda_d)$ that satisfy:
$$
\frac{x_1^2}{a_1-\lambda_i}+\dots+\frac{x_d^2}{a_d-\lambda_i}=1,\quad
1\le i\le d.
$$
\end{definition}

\begin{theorem}[Chasles]\label{th:Chasles}
Any line in $\mathbf E^d$ is tangent to exactly $d-1$ quadrics from
a given confocal family. Tangent hyper-planes to these quadrics,
constructed at the points of tangency with the line, are orthogonal
to each other.
\end{theorem}

\begin{definition}
We say that two lines $\ell_1$, $\ell_2$ intersecting in point
$P\in\mathcal Q$ {\em satisfy the billiard reflection law} on
quadric $\mathcal Q$ if they are coplanar with the line normal to
$\mathcal Q$ at $P$, determine the congruent angles with this line.

{\em A billiard trajectory within $\mathcal Q$} is a polygonal line
with vertices on $\mathcal Q$, such that each pair of its
consecutive segments satisfies the billiard reflection law on
$\mathcal Q$ and lies on the same side of the tangent hyper-plane to
$\mathcal Q$ at their common endpoint.
\end{definition}

The billiard reflection on a quadric is tightly connected with the
corresponding confocal family.

\begin{theorem}
Two lines that satisfy the reflection law on a quadric $\mathcal Q$
in $\mathbf E^d$ are tangent to the same $d-1$ quadrics confocal
with $\mathcal Q$.
\end{theorem}

From the previous theorem, it follows that all segment of a billiard
trajectory within $\mathcal Q$ are tangent to the same $d-1$
quadrics confocal with $\mathcal Q$. We call these $d-1$ quadrics
{\em caustics} of the given trajectory.

\begin{theorem} [Generalized Poncelet Theorem]
Consider a closed billiard trajectory within quadric $\mathcal Q$ in
$\mathbf E^d$. Then all other billiard trajectories within $\mathcal
Q$, that share the same $d-1$ caustics, are also closed. Moreover,
all these closed trajectories have the same number of vertices.
\end{theorem}

This theorem, concerning the plane case ($d=2$), is due to
Jean-Victor Poncelet and dates from the beginning of the XIXth
century \cite{Pon}. Its generalization to the three-dimensional case
is proved by Darboux \cite{Dar3} at the end of the same century. The
generalization to an arbitrary dimension is obtained in \cite{CCS}
at the end of XXth century.

\smallskip

Having the Generalized Poncelet Theorem in mind, it is natural to
ask about an explicit condition on a quadric and a set of confocal
caustics, such that the corresponding billiard trajectories are
closed. This condition, for $d=2$, was given by Cayley \cite{Cay}.
The condition for an arbitrary $d$ is obtained by the authors of
this paper in \cite{DR1,DR2}.

\begin{theorem} [Generalized Cayley Condition]
The condition on a billiard trajectory inside ellipsoid $\mathcal
Q_0$ in $\mathbf E^d$, with nondegenerate caustics $\mathcal
Q_{\alpha_1}$, \dots, $\mathcal Q_{\alpha_{d-1}}$ from the family
(\ref{eq:confocal.family}), to be perodic with period $n\ge d$ is:
\begin{equation*}
\rank\left(\begin{array}{cccc}
 B_{n+1} & B_n & \dots &  B_{d+1} \\
 B_{n+2} & B_{n+1} & \dots & B_{d+2} \\
 \dots \\
 \dots \\
 B_{2n-1} & B_{2n-2} & \dots & B_{n+d-1}
\end{array}\right)<n-d+1,
\end{equation*}
where
$\sqrt{(x-a_1)\dots(x-a_d)(x-\alpha_1)(x-\alpha_{d-1})}=B_0+B_1x+B_2x^2+\dots$
and all $a_1,\dots,a_d$ are distinct and positive.
\end{theorem}

\subsection{Poncelet Theorem in Projective Space over an Arbitrary
Field}\label{sekcija:proj}

In this section, we are going to present the generalization of
notions of confocal families of quadrics and related billiard
trajectories to the projective space over an arbitrary field of
characteristic not equal to $2$. Also, we give a proof of the Full
Poncelet theorem in the three-dimensional projective space, based on
\cite{CCS,CS,Ber}.

\smallskip

We begin with the definition of confocal family of quadrics in the
projective space, as it was done in \cite{CS}.

\begin{definition}\label{def:P.confocal}
A family of quadrics in the projective space is {\em confocal} if
for any hyper-plane the set of poles with respect to these quadrics
is a line.
\end{definition}

It can be shown that the dual to a confocal family is a pencil of
quadrics. In contrast to the Euclidean case, notice that it follows
from there that a confocal family in the projective space is
determined by two quadrics.

On the other hand, if we introduce the metrics in a convenient
manner in the real projective space $\mathbf {RP}^d$, the affine
part of a given confocal family in $\mathbf {RP}^d$ will satisfy
Definition \ref{def:E.confocal}. In this sense, Definitions
\ref{def:E.confocal} and \ref{def:P.confocal} are equivalent.

\smallskip

As well as in the Euclidean case, any line in $\mathbf P^d$ is also
tangent to $d-1$ quadrics from a given confocal family. We will
refer this statement by {\em Chasles theorem}, similarly as in the
Euclidean case, although the second part of Theorem \ref{th:Chasles}
has no sense without metrics.

\smallskip

Now, we are able to define the billiard reflection in a pure
projective manner, without metrics.

\smallskip

Let $\mathcal Q_1$ and $\mathcal Q_2$ be two quadrics that meet
transversely. Denote by $u$ the tangent plane to $\mathcal Q_1$ at
point $x$ and by $z$ the pole of $u$ with respect to $\mathcal Q_2$.
Suppose lines $\ell_1$ and $\ell_2$ intersect at $x$, and the plane
containing these two lines meet $u$ along $\ell$.

\begin{definition}
If lines $\ell_1, \ell_2, xz, \ell$ are coplanar and harmonically
conjugated, we say that rays $\ell_1$ and $\ell_2$ {\em obey the
reflection law} at the point $x$ of the quadric $\mathcal Q_1$ with
respect to the confocal family which contains $\mathcal Q_1$ and
$\mathcal Q_2$.
\end{definition}

If we introduce a coordinate system in which quadrics $\mathcal Q_1$
and $\mathcal Q_2$ are confocal in the Euclidean sense, reflection
defined in this way is same as the standard one.

\begin{theorem}[One Reflection Theorem]
Suppose rays $\ell_1$ and $\ell_2$ obey the reflection law at $x$ of
$\mathcal Q_1$ with respect to the confocal system determined by
quadrics $\mathcal Q_1$ and $\mathcal Q_2$. Let $\ell_1$ intersects
$\mathcal Q_2$ at $y_1'$ and $y_1$, $u$ is a tangent plane to
$\mathcal Q_1$ at $x$, and $z$ its pole with respect to $\mathcal
Q_2$. Then lines $y_1'z$ and $y_1z$ respectively contain
intersecting points $y_2'$ and $y_2$ of ray $\ell_2$ with $\mathcal
Q_2$. Converse is also true.
\end{theorem}

\begin{corollary} Let rays $\ell_1$ and $\ell_2$ obey the
reflection law of $\mathcal Q_1$ with respect to the confocal
system determined by quadrics $\mathcal Q_1$ and $\mathcal Q_2$.
Then $\ell_1$ is tangent to $\mathcal Q_2$ if and only if is
tangent $\ell_2$ to $\mathcal Q_2$; $\ell_1$ intersects $\mathcal
Q_2$ at two points if and only if $\ell_2$ intersects $\mathcal
Q_2$ at two points.
\end{corollary}

Next assertion is crucial for proof of the Poncelet theorem.

\begin{theorem} [Double Reflection Theorem] \label{th:DRT}
Suppose that $\mathcal Q_1$, $\mathcal Q_2$ are given quadrics and
$\ell_1$ line intersecting $\mathcal Q_1$ at the point $x_1$ and
$\mathcal Q_2$ at $y_1$. Let $u_1$, $v_1$ be tangent planes to
$\mathcal Q_1$, $\mathcal Q_2$ at points $x_1$, $y_1$ respectively,
and $z_1$, $w_1$ their poles with respect to $\mathcal Q_2$ and
$\mathcal Q_1$. Denote by $x_2$ second intersecting point of the
line $w_1x_1$ with $\mathcal Q_1$, by $y_2$ intersection of $y_1z_1$
with $\mathcal Q_2$ and by $\ell_2$, $\ell_1'$, $\ell_2'$ lines
$x_1y_2$, $y_1x_2$, $x_2y_2$. Then pairs $\ell_1,\ell_2$;
$\ell_1,\ell_1'$; $\ell_2,\ell_2'$; $\ell_1',\ell_2'$ obey the
reflection law at points $x_1$ (of $\mathcal Q_1$), $y_1$ (of
$\mathcal Q_2$), $y_2$ (of $\mathcal Q_2$), $x_2$ (of $\mathcal
Q_1$) respectively.
\end{theorem}

\begin{corollary}
If the line $\ell_1$ is tangent to a quadric $\mathcal Q$ confocal
with $\mathcal Q_1$ and $\mathcal Q_2$, then rays $\ell_2$,
$\ell_1'$, $\ell_2'$ also touch $\mathcal Q$.
\end{corollary}

\medskip

Using Double reflection theorem, we can prove the Full Poncelet
theorem in the three-dimensional space. Let $\mathcal Q_1, \dots,
\mathcal Q_m, \mathcal Q$ be confocal quadrics and $\ell_1, \dots,
\ell_m$, ($\ell_{m+1} = \ell_1$) lines such that pairs $\ell_i,
\ell_{i+1}$ obey the reflection law at point $x_i$ of $\mathcal
Q_i$. Let $u_i$ be a plane tangent to $\mathcal Q_i$ at $x_i$, $z_i$
its pole with respect to $\mathcal Q$ and $z_i \not\in\mathcal Q$.
If line $\ell_1$ intersects $\mathcal Q$ at $y_1, y_1'$, then, by
One reflection theorem, the second intersection point $y_2$ of line
$y_1z_1$ with quadrics $\mathcal Q$ belongs to $\ell_2$. Similarly,
having point $y_i$, we construct the point $y_{i+1} \in \ell_{i+1}$.
It follows that $y_{m+1} = y_1$ or $y_{m+1} = y_1'$. Suppose
$y_{m+1} = y_1$. It can be proved that, if, for a given polygon $x_1
\dots x_m$ and quadric $\mathcal Q$, $y_{m+1}=y_1$ holds, then
$y_{m+1}=y_1$ for any surface $\mathcal Q$ from the confocal family.

\smallskip

Suppose $\ell_i'$ are rays $\ell_i$ reflected of $\mathcal Q$ at
points $y_i$. By Double reflection theorem, $\ell_i'$ and
$\ell_{i+1}'$ meet at point $x_i' \in\mathcal Q_i$ and obey the
reflection law. In this way, we obtained a new polygon $x_1' \dots
x_m'$. If $\ell_1$ touches quadrics $\mathcal Q',\mathcal Q''$
confocal with $\{\mathcal Q_1,\mathcal Q \}$, then all sides of both
polygons touch them. In this way, a two-parameter family of Poncelet
polygons is obtained.

\subsection{Virtual Billiard Trajectories}\label{virtual}

Apart from the real motion of the billiard particle in $\mathbf
E^d$, it is of interest to consider {\it virtual reflections}. These
reflections were discussed by Darboux in \cite{Dar3} (see Chapter
XIV of Book IV in Volume 2). In this section, we review some recent
results from \cite{DR4}.

\smallskip

\begin{figure}[h]
\centering
\begin{minipage}[b]{0.6\textwidth}
\centering
\includegraphics[width=6cm,height=5cm]{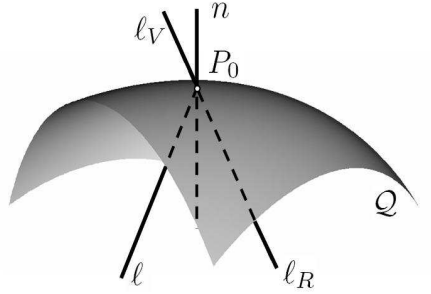}
\caption{Real and virtual
reflection}\label{fig:virtual.real.reflection}
\end{minipage}
\end{figure}

Formally, in the Euclidean space, we can define the {\it virtual
reflection} at the quadric $\mathcal Q$ as a map of a ray $\ell$
with the endpoint $P_0$ ($P_0\in\mathcal Q$) to the ray
complementary to the one obtained from $\ell$ by the real reflection
from $\mathcal Q$ at the point $P_0$. On Figure
\ref{fig:virtual.real.reflection}, ray $\ell_R$ is obtained from
$\ell$ by the real reflection on the quadric surface $\mathcal Q$ at
point $P_0$. Ray $\ell_V$ is obtained from $\ell$ by the virtual
reflection. Line $n$ is the normal to $\mathcal Q$ at $P_0$.

\smallskip

Let us remark that, in the case of real reflections, exactly one
elliptic coordinate, the one corresponding to the quadric $\mathcal
Q$, has a local extreme value at the point of reflection. On the
other hand, on a virtual reflected ray, this coordinate is the only
one not having a local extreme value at the point of reflection. In
the $2$-dimensional case, the virtual reflection can easily be
described as the real reflection from the other confocal conic
passing through the point $P_0$. In higher dimensional cases, the
virtual reflection can be regarded as the real reflection from the
line normal to $\mathcal Q$ at $P_0$ (see Figure
\ref{fig:virtual.real.reflection}).

\smallskip

The notions of real and virtual reflection cannot be
straightforwardly extended to the projective space, since there we
essentially use that the field of real numbers is naturally ordered.
Nevertheless, it turns out that it is possible to define a certain
configuration connected with real and virtual reflection, such that
its properties remain in the projective case, too.

\smallskip

Let points  $X_1, X_2$; $Y_1, Y_2$ belong to quadrics $\mathcal
Q_1$, $\mathcal Q_2$ in $\mathbf P^d$.

\begin{definition}\label{def:VRC}
We will say that the quadruple of points $X_1, X_2, Y_1, Y_2$
constitutes a {\em virtual reflection configuration (VRC)} if pairs
of lines $X_1 Y_1$, $X_1 Y_2$; $X_2 Y_1$, $X_2 Y_2$; $X_1 Y_1$, $X_2
Y_1$; $X_1 Y_2$, $X_2 Y_2$ satisfy the reflection law at points
$X_1$, $X_2$ of $\mathcal Q_1$ and $Y_1$, $Y_2$ of $\mathcal Q_2$
respectively, with respect to the confocal system determined by
$\mathcal Q_1$ and $\mathcal Q_2$.

If, additionally, the tangent planes to $\mathcal Q_1, \mathcal Q_2$
at $X_1, X_2$; $Y_1, Y_2$ belong to a pencil, we say that these
points constitute a {\em double reflection configuration (DRC)} (see
Figure \ref{fig:virtual_reflection}).
\begin{figure}[h]
\centering
\begin{minipage}[b]{0.6\textwidth}
\centering
\includegraphics[width=5cm,height=4cm]{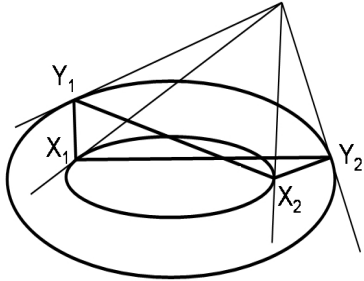}
\caption{Double reflection configuration}
\end{minipage}\label{fig:virtual_reflection}
\end{figure}
\end{definition}

Now, the Darboux's statement can be generalized and proved as
follows:

\begin{theorem}\label{th:virt.refl}
Let $\mathcal Q_1$, $\mathcal Q_2$ be two quadrics in the projective
space $\mathbf P^d$, $X_1$, $X_2$ points on $\mathcal Q_1$ and
$Y_1$, $Y_2$ on $\mathcal Q_2$. If the tangent hyperplanes at these
points to the quadrics belong to a pencil, then $X_1, X_2, Y_1, Y_2$
constitute a virtual reflection configuration.

Furthermore, suppose that the projective space is defined over the
field of reals. Introduce a coordinate system, such that $\mathcal
Q_1$, $\mathcal Q_2$ become confocal ellipsoids in the Euclidean
space. If $\mathcal Q_1$ is placed inside $\mathcal Q_2$, then the
sides of the quadrilateral $X_1Y_1X_2Y_2$ obey the real reflection
from $\mathcal Q_2$ and the virtual reflection from $\mathcal Q_2$.
\end{theorem}

We are going to conclude this section with the statement converse to
the previous theorem.

\begin{proposition}
Let pairs of points $X_1$, $X_2$ and $Y_1$, $Y_2$ belong to confocal
ellipsoids $\mathcal Q_1$ and $\mathcal Q_2$ in Euclidean space
$\mathbf E^d$, and let $\alpha_1$, $\alpha_2$, $\beta_1$, $\beta_2$
be the corresponding tangent planes. If a quadruple $X_1, X_2, Y_1,
Y_2$ is a virtual reflection configuration, then planes $\alpha_1$,
$\alpha_2$, $\beta_1$, $\beta_2$ belong to a pencil.
\end{proposition}

\subsection{Generalized Cayley Curve}\label{gen.lebeg}

\begin{definition}\label{def:gen.Cayley}
Let $\ell$ be a line not contained in any quadric of the given
confocal family in the projective space $\mathbf P^d$. The {\em
generalized Cayley curve} $\mathcal C_{\ell}$ is the variety of
hyperplanes tangent to quadrics of the confocal family at the points
of $\ell$.
\end{definition}

This curve is naturally embedded in the dual space
$\mathbf{P}^{d\,*}$.

\smallskip

\begin{figure}[h]
\centering
\begin{minipage}{0.6\textwidth}
\centering
\includegraphics[width=7cm,height=4cm]{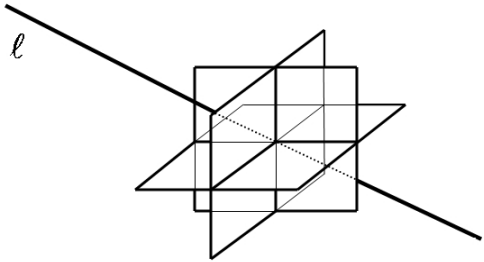}
\caption{Three points of the generalized Cayley curve in dimension
3}\label{fig:gen.cayley.curve}
\end{minipage}
\end{figure}

On Figure \ref{fig:gen.cayley.curve}, we see the planes which
correspond to one point of the line $\ell$ in the $3$-dimensional
space.

\begin{proposition}
The generalized Cayley curve is a hyperelliptic curve of genus
$g=d-1$, for $d\ge3$. Its natural realization in $\mathbf P^{d\,*}$
is of degree $2d-1$.
\end{proposition}

The natural involution $\tau_{\ell}$ on the generalized Cayley's
curve $\mathcal C_{\ell}$ maps to each other the tangent planes at
the points of intersection of $\ell$ with any quadric of the
confocal family. It is easy to see that the fixed points of this
involution are hyperplanes corresponding to the $d-1$ quadrics that
are touching $\ell$ and to $d+1$ degenerate quadrics of the confocal
family.

\smallskip

Now, let us make a few remarks on the real case.

\begin{lemma}\label{lema:polinom}
Suppose a line $\ell$ is tangent to quadrics
 $\mathcal Q_{\alpha_1},\dots,\mathcal Q_{\alpha_{d-1}}$
from the family {\rm (\ref{eq:confocal.family})}. Then Jacobian
coordinates $(\lambda_1,\dots, \lambda_d)$ of any point on $\ell$
satisfy the inequalities $\mathcal P(\lambda_s)\ge 0$,
$s=1,\dots,d$, where
\begin{equation}\label{eq:polinom}
\mathcal P(x)=(a_1-x)\dots(a_d-x)(\alpha_1-x)\dots(\alpha_{d-1}-x).
\end{equation}
\end{lemma}

\begin{proof}
Let $x$ be a point of $\ell$, $(\lambda_1,\dots,\lambda_d)$ its
Jacobian coordinates, and $y$ a vector parallel to $\ell$. The
equation $Q_{\lambda}(x+ty)=1$ is quadratic with respect to $t$. Its
discriminant is:
$$
\Phi_{\lambda}(x,y) =
Q_{\lambda}(x,y)^2-Q_{\lambda}(y)\bigl(Q_{\lambda}(x)-1\bigr),
$$
where
$$
Q_{\lambda}(x,y) = \frac{x_1y_1}{a_1-\lambda}+\dots +
\frac{x_dy_d}{a_d-\lambda}.
$$
By \cite{Mo},
$$
\Phi_{\lambda}(x,y)=\frac
{(\alpha_1-\lambda)\dots(\alpha_{d-1}-\lambda)}
{(a_1-\lambda)\dots(a_d-\lambda)}.
$$
For each of the coordinates $\lambda=\lambda_s$, ($1\le s\le d$),
the quadratic equation has a solution $t=0$; thus, the corresponding
discriminants are non-negative. This is obviously equivalent to
$\mathcal P(\lambda_s)\ge0$.
\end{proof}

The equation of the generalized Cayley curve corresponding to a
confocal family of the form (\ref{eq:confocal.family}) can be
written as:
$$
y^2=\mathcal P(x).
$$
It is important to note that the constants $\alpha_1$, \dots,
$\alpha_{d-1}$, corresponding to the quadrics that are touching
$\ell$ cannot take arbitrary values. More precisely, following
\cite{Au,Kn}, we can state:

\begin{proposition}
There exists a line in $\mathbf E^d$ that is tangent to $d-1$
distinct non-degenerate quadrics $\mathcal Q_{\alpha_1}$, \dots,
$\mathcal Q_{\alpha_{d-1}}$ from the family
(\ref{eq:confocal.family}) if and only if the set
$\{a_1,\dots,a_d,\alpha_1,\dots,\alpha_{d-1}\}$ can be ordered as
$b_1<b_2<\dots<b_{2d-1}$, such that $\alpha_j\in\{b_{2j-1},b_{2j}\}$
$(1\le j\le d-1)$.
\end{proposition}

It was observed in \cite{DR4} that the generalized Cayley curve is
isomorphic to the Veselov-Moser isospectral curve.

\smallskip

This curve is also naturally isomorphic to the curves studied by
Kn\"orrer, Donagi and Reid (RDK curves). We begin next section with
the construction of this isomorphism, in order to establish the
relationship between billiard law and algebraic structure on the
Jacobian $\Jac(\mathcal C_{\ell})$.

\section{Billiard Law and Algebraic Structure on the Abelian Variety $\mathcal
A_{\ell}$}\label{sekcija:bilijarskaalgebra}

The aim of this section is to construct in $\mathcal A_{\ell}$ an
algebraic structure that is naturally connected with the billiard
motion. In Section \ref{sekcija:morfizam}, we show that the
generalized Cayley's curve is isomorphic to the RDK curve. Then, in
order to get better understanding and intuition, we are going first,
in Section \ref{sec:genus2} to describe in detail the billiard
algebra, and prove some of its nice geometrical properties, for the
case of dimension $3$, i.e.\ when the corresponding curve is of
genus $2$. The general construction is given in the Section
\ref{sec:general}.

\subsection{Morphism between Generalized Cayley's Curve and RDK
Curve}\label{sekcija:morfizam}

Now we are going to establish the connection between generalized
Cayley's curve defined above and the curves studied by Kn\"orrer,
Donagi, Reid and to trace out the relationship between billiard
constructions and the algebraic structure of the corresponding
Abelian varieties.

\smallskip

Suppose a line $\ell$ in $\mathbf E^d$ is tangent to quadrics
 $\mathcal Q_{\alpha_1},\dots,\mathcal Q_{\alpha_{d-1}}$
from the confocal family (\ref{eq:confocal.family}). Denote by
$\mathcal A_{\ell}$ the family of all lines which are tangent to the
same $d-1$ quadrics. Note that according to the corollary of the One
Reflection Theorem, the set $\mathcal A_{\ell}$ is invariant to the
billiard reflection on any of the confocal quadrics.

\smallskip

We begin with the next simple observation.

\begin{lemma}\label{lema:veza}
Let the lines $\ell$ and $\ell'$ obey the reflection law at the
point $z$ of a quadric $\mathcal Q$ and suppose they are tangent to
a confocal quadric $\mathcal Q_1$ at the points $z_1$ and $z_2$.
Then the intersection of the tangent spaces $T_{z_1}\mathcal Q_1\cap
T_{z_2}\mathcal Q_1$ is contained in the tangent space
$T_{z}\mathcal Q$.
\end{lemma}

\begin{proof}
It follows from the One Reflection Theorem: since the poles $z_1$,
$z_2$ and $w$ of the planes $T_{z_1}\mathcal Q_1$,
$T_{z_2}\mathcal Q_1$ and $T_z\mathcal Q$ with respect to the
quadric $\mathcal Q_1$ are colinear, the planes belong to a
pencil.
\end{proof}

Following \cite{Kn}, together with $d-1$ affine confocal quadrics
$\mathcal Q_{\alpha_1},\dots,\mathcal Q_{\alpha_{d-1}}$, one can
consider their projective closures $\mathcal
Q^p_{\alpha_1},\dots,\mathcal Q^p_{\alpha_{d-1}}$ and the
intersection $V$ of two quadrics in $\mathbf P^{2d-1}$:
\begin{equation}\label{eq:kvadrika1}
x_1^2+\dots+x_d^2-y_1^2-\dots - y_{d-1}^2=0,
\end{equation}
\begin{equation}\label{eq:kvadrika2}
 a_1x_1^2+\dots +a_dx_d^2 - \alpha_1y_1^2 - \dots -\alpha_{d-1}y_{d-1}^2=x_0^2.
\end{equation}

Denote by $F=F(V)$ the set of all $(d-2)$-dimensional linear
subspaces of $V$. For a given $L\in F$, denote by $F_L$ the closure
in $F$ of the set
 $\{\, L'\in F \mid \dim L\cap L' = d-3\, \}$.
It was shown in \cite{Re} that $F_L$ is a nonsingular hyperelliptic
curve of genus $d-1$. Note that for $d=3$, i.e.\ when the curve
$\mathcal C_{\ell}$ is of the genus $2$, an isomorphism between
$F(V)$ and the Jacobian of the hyperelliptic curve was established
in \cite{NR}.

\smallskip

The projection
$$
\pi'\ :\ \mathbf P^{2d-1}\setminus\{(x,y)|x=0\}\to\mathbf P^d, \quad
\pi'(x,y)=x,
$$
maps $L\in F(V)$ to a subspace $\pi'(L)\subset\mathbf P^d$ of the
codimension $2$. $\pi'(L)$ is tangent to the quadrics $\mathcal
Q^{p*}_{\alpha_1},\dots,\mathcal Q^{p*}_{\alpha_{d-1}}$ that are
dual to $\mathcal Q^p_{\alpha_1},\dots,\mathcal Q^p_{\alpha_{d-1}}$.

\smallskip

Thus, the space dual to $\pi'(L)$, denoted by $\pi^*(L)$, is a line
tangent to the quadrics $\mathcal Q^p_{\alpha_1},\dots,\mathcal
Q^p_{\alpha_{d-1}}$.

\smallskip

We can reinterpret the generalized Cayley's curve $\mathcal
C_{\ell}$, which is a family of tangent hyperplanes, as a set of
lines from $\mathcal A_{\ell}$ which intersect $\ell$. Namely, for
almost every tangent hyperplane there is a unique line $\ell'$,
obtained from $\ell$ by the billiard reflection. Having this
identification in mind, it is easy to prove the following

\begin{corollary}
There is a birational morphism between the generalized Cayley's
curve $\mathcal C_{\ell}$ and Reid-Donagi-Kn\"orrer's curve $F_L$,
with $L=\pi^{*-1}(\ell)$, defined by
$$
j: \ell'\mapsto L',\quad L'=\pi^{*-1}(\ell'),
$$
where $\ell'$ is a line obtained from $\ell$ by the billiard
reflection on a confocal quadric.
\end{corollary}

\begin{proof}
It follows from the previous lemmata and Lemma 4.1 and Corollary 4.2
from \cite{Kn}.
\end{proof}

Thus, Lemma \ref{lema:veza} gives a link between the dynamics of
ellipsoidal billiards and algebraic structure of certain Abelian
varieties. This link provides a two way interaction: to apply
algebraic methods in the study of the billiard motion, but also vice
versa, to use billiard constructions in order to get more effective,
more constructive and more observable understanding of the algebraic
structure.

\smallskip

In the following section, we are going to use this link in
constructing {\em the billiard algebra}, which is a group structure
in $\mathcal A_{\ell}$.

\subsection{Genus 2 Case}\label{sec:genus2}

Before we proceed in general case, we want to emphasize the billiard
constructions involved in the first nontrivial case, of genus two.

\subsubsection*{Leading Principle, Definition and First Properties of the
Operation}

We formulate {\em The Leading Principle}:

\smallskip

{\em The sum of the lines in any virtual reflection configuration is
equal to zero if the four tangent planes at the points of reflection
belong to a pencil.}

\smallskip

Recall that, by Definition \ref{def:VRC}, such a configuration of
four lines in a VRC, with tangent planes in a pencil, we are going
to call {\em a Double Reflection Configuration (DRC)}.

\medskip

{\bf Neutral element}

\smallskip

Let us fix a line $\mathcal O\in\mathcal A_{\ell}$.

\medskip

{\bf Opposite element}

\smallskip

First, we define $-\mathcal O:=\mathcal O$.

\smallskip

For a given line $x\in\mathcal C_{\mathcal O}$, define
$$
-x:=\tau_{\mathcal O}(x),
$$
where $\tau_{\mathcal O}$ is the hyperelliptic involution of the
curve $\mathcal C_{\mathcal O}$.

\begin{proposition}\label{prop:saglasnost-}
For any $x\in\mathcal C_{\mathcal O}$, both lines $x$ and $-x$ are
obtained from $\mathcal O$ by the reflection on the same quadric
$\mathcal Q_x$ from the confocal family.

\smallskip

Moreover, let $\pi$ be the unique plane orthogonal to $\mathcal O$
from the pencil determined by the tangent planes to $\mathcal Q_x$
at the intersection points with $\mathcal O$, and $\mathcal
Q_{\mathcal O}$ the unique quadric from the confocal family such
that $\pi$ is tangent to it. Then the intersection point of $\pi$
with $\mathcal O$ belongs to $\mathcal Q_{\mathcal O}$.
\end{proposition}

\begin{proof}
Follows from the degenerate case of Double Reflection Theorem
applied on quadrics $\mathcal Q_x$ and $\mathcal Q_{\mathcal O}$
with $\mathcal O = l_1 = l_2$. See Figure \ref{fig:saglasnost-}.
\begin{figure}[h]
\centering
\begin{minipage}{0.6\textwidth}
\centering
\includegraphics[width=6cm,height=4cm]{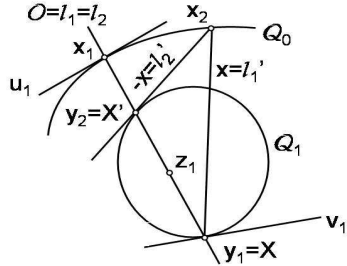}
\caption{Proposition \ref{prop:saglasnost-}}\label{fig:saglasnost-}
\end{minipage}
\end{figure}
\end{proof}

\begin{corollary}
Plane $\pi$ and quadric $\mathcal Q_{\mathcal O}$ do not depend on
$x$. Line $\mathcal O$ reflects to itself on $\mathcal Q_{\mathcal
O}$.
\end{corollary}

Suppose now that $x\in \mathcal A_{\ell}$, but $x$ doesn't belong to
$\mathcal C_{\mathcal O}$. Thus $x$ does not intersect $\mathcal O$
and the two lines generate a linear projective space. This space
intersects $\mathcal A_{\ell}$ along the divisor
$$
\mathcal O + x + p + q.
$$
(See, for example, \cite{Tyu}, \cite{GH3}.)

\smallskip

According to the Double Reflection Theorem, one can see that the
lines $\mathcal O$, $x$, $p$, $q$ form a double reflection
configuration. We define $-x$ such that it forms a double reflection
configuration with $\mathcal O$, $-p$, $-q$.

\smallskip

From the definition, it is immediately seen that $-(-x)=x$ for every
$x\in\mathcal A_{\ell}$.

\smallskip

The following property is a consequence of Proposition
\ref{prop:saglasnost-} and the Double Reflection Theorem.

\begin{proposition}\label{prop:x.-x}
Lines $x$ and $-x$ intersect each other and they satisfy the
reflection law on $\mathcal Q_{\mathcal O}$.
\end{proposition}

The following example is an illustration for the construction we
have just made.

\begin{example}
Take the line $\mathcal O$ to be orthogonal to one of the coordinate
hyperplanes. Then $\mathcal Q_{\mathcal O}=\pi$ coincide with this
hyperplane and, additionally, among caustics $\mathcal Q_1$, \dots,
$\mathcal Q_{d-1}$ there cannot be two quadrics of the same type.

\smallskip

According to \cite{Au}, exactly the case of all quadrics of
different type corresponds to the case where $\mathcal A_{\ell}$
consists of only one real connected component, isomorphic to the
connected component of zero of $\Jac(\mathcal C_{\mathcal
O})(\mathbf R)$.

\smallskip

In this case, for any line $x\in\mathcal A_{\ell}$, the opposite
element $-x$ may be defined as the line symmetric to $x$ with
respect to plane $\mathcal Q_{\mathcal O}$.
\end{example}

\medskip

{\bf Addition}

\smallskip

We are going to define operation
$$
+\ :\ \mathcal A_{\ell}\times \mathcal A_{\ell} \to \mathcal
A_{\ell}.
$$

\smallskip

Define $\mathcal O+x = x+\mathcal O=x$, for all $x\in\mathcal
A_{\ell}$.

\smallskip

For $s_1,s_2\in\mathcal C_{\mathcal O}$, define $s_1+s_2$ as the
line that forms a double reflection configuration with $-s_1, -s_2,
\mathcal O$. Obviously, $s_1+s_2=s_2+s_1$.

\smallskip

Notice that $-s_1,-s_2$ are unique lines from $\mathcal A_{\ell}$
that intersect both $s_1+s_2$ and $\mathcal O$ (see \cite{Tyu}).
Thus, we have:

\begin{lemma}\label{lema:rastavljanje}
Each line $x\in\mathcal A_{\ell}\setminus\mathcal C_{\mathcal O}$
can be in the unique way represented as the sum of two lines that
intersect $\mathcal O$.
\end{lemma}

\begin{figure}[h]
\centering
\begin{minipage}{0.6\textwidth}
\centering
\includegraphics[width=6cm,height=4.5cm]{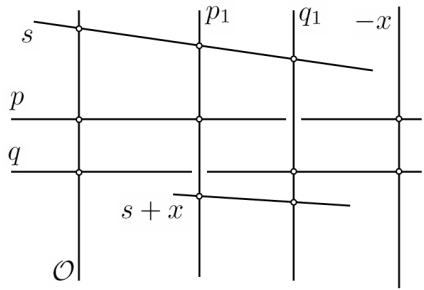}
\caption{Partial operation}\label{fig:partial1}
\end{minipage}
\end{figure}

Now, suppose $s\in\mathcal C_{\mathcal O}$, $x\in\mathcal
A_{\ell}\setminus\mathcal C_{\mathcal O}$ (see Figure
\ref{fig:partial1}).
As already explained, take lines $p,q\in\mathcal C_{\mathcal O}$
such that $x=p+q$. Construct
$$
p_1= (-s) + (-p),\quad q_1=(-s)+(-q)
$$
as above, since both pairs $s, p$ and $s, q$ intersect $\mathcal O$.
Now both $p_1$ and $q_1$ intersect $s$. Thus, the three lines belong
to a DRC with the fourth line $z$. We put by definition
$$
s+x=x+s=z.
$$

The following lemma gives a very important property of the
operation.

\begin{lemma}\label{lema:saglasnost}
Let $s,x$ be lines in $\mathcal C_{\mathcal O}$, $\mathcal A_{\ell}$
respectively and $\mathcal Q_s$ the quadric from the confocal family
such that $s$ and $\mathcal O$ reflect to each other on it. Then the
lines $s+x$ and $-x$ intersect each other, their intersection point
belongs to quadric $\mathcal Q_s$, and these two lines satisfy the
billiard reflection law on $\mathcal Q_s$.
\begin{figure}[h]
\centering
\begin{minipage}{0.6\textwidth}
\centering
\includegraphics[width=6cm,height=4.5cm]{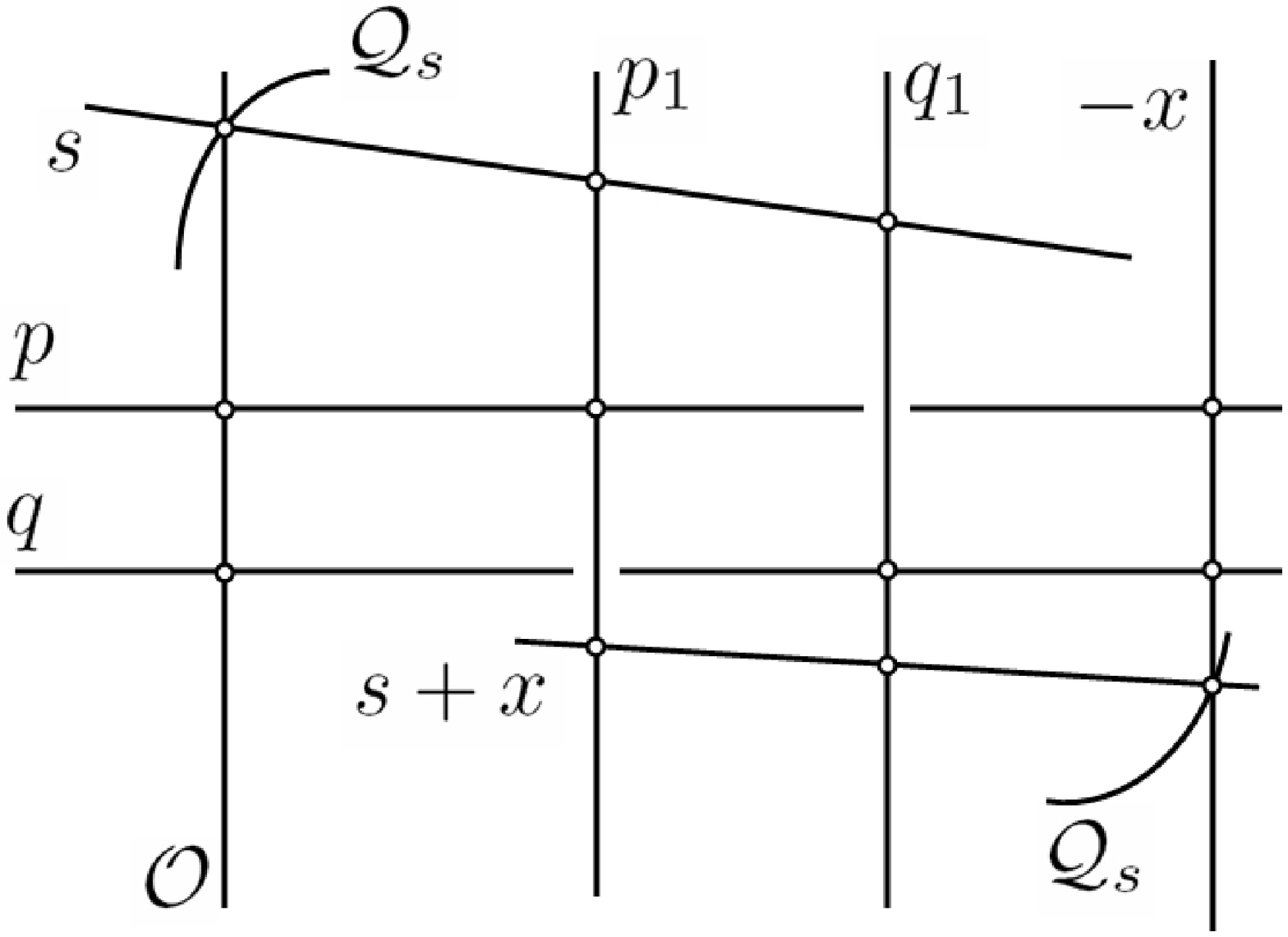}
\caption{Lemma \ref{lema:saglasnost}}\label{fig:saglasnost}
\end{minipage}
\end{figure}
\end{lemma}

\begin{proof}
Follows from \cite{Don}, \cite{Kn} and Lemma \ref{lema:veza}.
\end{proof}

\medskip

Lemma \ref{lema:saglasnost} is going to play an important role in
proving basic properties of the operation:

\begin{lemma}\label{lema:asoc}
For all $p,q,s\in\mathcal C_{\mathcal O}$, the associative law
holds: $p+(q+s)=(p+q)+s$.
\end{lemma}

\begin{proof}
Denoting, like on Figure \ref{fig:saglasnost}, by $-x$, $p_1$, $q_1$
lines forming DRCs with triplets $(\mathcal O,p,q)$, $(\mathcal
O,p,s)$, $(\mathcal O,q,s)$, and applying Lemma
\ref{lema:saglasnost}, we see that $(p+q)+s$ is the unique line
forming a DRC with triplets $(p_1,q_1,s)$, $(p,p_1,-x)$,
$(q,q_1,-x)$. The same holds for line $p+(q+s)$, thus the two lines
coincide.
\end{proof}

\medskip

\begin{lemma}\label{lema:-(p+q+s)=(-p)+(-q)+(-s)}
Let $p,q,s\in\mathcal C_{\mathcal O}$. Then
$-(p+q+s)=(-p)+(-q)+(-s)$.
\end{lemma}

\begin{proof}
It is straightforward to prove that $p+q+s\in\mathcal C_{\mathcal
O}$ if and only if two of the lines $p,q,s$ are inverse to each
other. In this case, the equality we need to prove immediately
follows.

\smallskip

Thus, suppose $p+q+s$ does not intersect $\mathcal O$ and that $a,b$
are lines forming a DRC with $p+q+s,\mathcal O$. By Lemma
\ref{lema:saglasnost}, we have that $p+q+s$ is the unique line
intersecting $(-p)+(-q)$, $(-p)+(-s)$ and $(-q)+(-s)$. Reflecting
$p,q,s,a,b$ on $\mathcal Q_{\mathcal O}$ and applying the Double
Reflection Theorem, we get that $-(p+q+s)$ intersects lines $p+q$,
$p+s$, $q+s$. Thus, it is equal to $(-p)+(-q)+(-s)$, by Lemma
\ref{lema:saglasnost}.
\end{proof}

\medskip

Suppose now two lines $x, y \in \mathcal A_{\ell}$ are given, none
of which is intersecting $\mathcal O$. We define their sum as
follows.

\smallskip

First, we represent $x$ as a sum of two lines intersecting $\mathcal
O$: $x=s_1+s_2$, $s_1,s_2\in\mathcal C_{\mathcal O}$. Then, we
define:
$$
x+y:=s_1+(s_2+y).
$$

\smallskip

We need to show that this definition is correct.

\begin{lemma}\label{lema:dobrodef}
Let $s_1,s_2\in\mathcal C_{\mathcal O}$, and $y\in\mathcal
A_{\ell}$. Then
$$
s_1 + (s_2 +y)=s_2 + (s_1+y).
$$
\end{lemma}

\begin{proof}
If $y\in\mathcal C_{\mathcal O}$, it is enough to apply Lemma
\ref{lema:asoc} and the commutativity property for the addition in
$\mathcal C_{\mathcal O}$.

\smallskip

If $y\in\mathcal A_{\ell}\setminus\mathcal C_{\mathcal O}$, then the
lines $-s_2-y$ and $y$ intersect with billiard reflection, as well
as lines $-s_1-y$ and $y$ do. Thus, there is a unique, fourth line
in the intersection of the space generated by $[-s_1-y,y, -s_2-y]$
and $\mathcal A_{\ell}$. At one hand, this line is equal to
$s_2-(-s_1-y)$, at the other it is equal to $s_1-(-s_2-y)$.
\end{proof}

\smallskip

Now, from Lemmata \ref{lema:asoc}--\ref{lema:dobrodef}, it follows
that we have constructed on $\mathcal A_{\ell}$ a commutative group
structure that is naturally connected with the billiard law.

\subsubsection*{Further Properties of the Operation}

We started this section with The Leading Principle, and used it as a
natural motivation for the construction of the algebra. Now, we are
going to show that this is justified, i.e.\ that the statement
formulated as The Leading Principle really holds in the group
structure.

\begin{theorem}\label{th:leading}
The sum of the lines in any double reflection configuration is equal
to zero.
\end{theorem}

\begin{proof}
Let $a,b,c,d$ be the lines of a DRC. Take that pairs $a,b$ and $c,d$
satisfy the reflection law on $\mathcal Q_1$, and pairs $b,c$, $d,a$
on $\mathcal Q_2$. Then, by Lemma \ref{lema:saglasnost}, we have:
\begin{equation}\label{eq:drc}
b=-a+s_1,\quad c=-b+s_2,\quad d=-c+\bar s_1,\quad a=-d+\bar s_2,
\end{equation}
where $s_i$, $\bar s_i$  are obtained from $\mathcal O$ by the
billiard reflection on $\mathcal Q_i$, $i=1,2$. Obviously, $\bar
s_i\in\{s_i,-s_i\}$.

\smallskip

From (\ref{eq:drc}), we get:
\begin{equation}\label{eq:s1s2}
a+b+c+d=s_1+\bar s_1=s_2+\bar s_2.
\end{equation}
Thus, we only need to check the case when $\bar s_1=s_1$ and $\bar
s_2=s_2$. Then we have: $s_2=s_1+s_1+(-s_2)$ and, from the
definition of the addition operation, it follows that lines
$\mathcal O,s_1,-2s_1, s_2$ constitute a closed billiard trajectory
with consecutive reflections on
 $\mathcal Q_1,\mathcal Q_1,\mathcal Q_2,\mathcal Q_2$.
On the other hand, $(-s_1)+(-s_2)$ is the unique line in $\mathcal
A_{\ell}$ that, besides $\mathcal O$, intersects both $s_1$ and
$s_2$. Thus,
$$
-2s_1=(-s_1)+(-s_2)\ \Rightarrow\ s_1=s_2\ \Rightarrow\ \mathcal
Q_1=\mathcal Q_2.
$$
This means that the double reflection configuration $a,b,c,d$ is
degenerated, i.e.\ two of the lines in the configuration coincide.
Say $a=c$, then both $b$ and $d$ are obtained from $a$ by the
billiard reflection on $\mathcal Q_1$. So, $b=-a+s_1$, $d=-a-s_1$
and the theorem follows.
\end{proof}

\medskip

Let us consider a billiard trajectory $\mathbf
t=(\ell_0,\ell_1,\dots,\ell_n)$ with the initial line
$\ell_0=\mathcal O$ and the reflections on quadrics $\mathcal
Q_1,\dots,\mathcal Q_n$ from the confocal family. Applying Double
Reflection Theorem to $\mathbf t$, we obtain different trajectories
with sharing initial segment $\mathcal O$ and final segment
$\ell_n$. The reflections on each such a trajectory are also on the
same set of quadrics $\mathcal Q_1,\dots,\mathcal Q_n$, but their
order may be changed. Let us denote by $\ell_1^{(1)}$, \dots,
$\ell_n^{(n)}$ all possible lines obtained after the first
reflection on all these trajectories. Then, in our algebra, the
final segment $\ell_n$ can be calculated from these ones:

\begin{proposition}\label{prop:elln}
$\ell_n=(-1)^{n+1}(\ell_1^{(1)}+\dots+\ell_n^{(n)}).$
\end{proposition}

\begin{proof}
Follows from Theorem \ref{th:leading}.
\end{proof}

\medskip

The following interesting property also may be proved from Lemmata
\ref{lema:asoc}--\ref{lema:dobrodef}.

\begin{proposition}
Let line $x\in\mathcal{A}_{\ell}$ be obtained from $\mathcal O$ by
consecutive reflections on three quadrics $\mathcal Q_1$, $\mathcal
Q_2$, $\mathcal Q_3$. Then the $6$ lines obtained from $\mathcal O$
by reflections of the these quadrics may be divided into $2$ groups
such that:

-- for each $i\in\{1,2,3\}$, the lines obtained from $\mathcal O$ by
the reflections on $\mathcal Q_i$ are in different groups;

-- all trajectories with three reflections on $\mathcal Q_1$,
$\mathcal Q_2$, $\mathcal Q_3$ starting with $\mathcal O$ and ending
with $x$ contain lines only from one of the groups.
\end{proposition}

We will finish this subsection with a very beautiful and non-trivial
theorem on confocal families of quadrics.

\begin{theorem}\label{th:zvezda}
Let $\mathcal F$ be a family of confocal quadrics in $\mathbf P^3$.
There exist configurations consisting of $12$ planes in $\mathbf
P^3$ with the following properties:
\begin{itemize}
\item[-]
The planes may be organized in $8$ triplets, such that each plane in
a triplet is tangent to a different quadric from $\mathcal F$ and
the three touching points are collinear. Every plane in the
configuration is a member of two triplets.
\item[-]
The planes may be organized in $6$ quadruplets, such that the planes
in each quadruplet belong to a pencil and they are tangent to two
different quadrics from $\mathcal F$. Every plane in the
configuration is a member of two quadruplets.
\end{itemize}
Moreover, such a configuration is determined by three planes tangent
to three different quadrics from $\mathcal F$, with collinear
touching points.
\end{theorem}

\begin{proof}
Denote by $\mathcal O$ the line containing the three touching points
and $p,q,s$ the lines obtained from $\mathcal O$ by the billiard
reflection on the given quadrics from $\mathcal F$. We construct
lines $p_1,q_1,-x,x+s$ as explained before Lemma
\ref{lema:saglasnost} (see Figure \ref{fig:saglasnost}). The planes
of the configuration are tangent to corresponding quadrics at points
of the intersection of the lines.

\smallskip

The configuration of the planes in the dual space $\mathbf P^{3*}$
is shown on Figure \ref{fig:zvezda}.
\begin{figure}[h]
\centering
\begin{minipage}{0.6\textwidth}
\centering
\includegraphics[width=6cm,height=4.5cm]{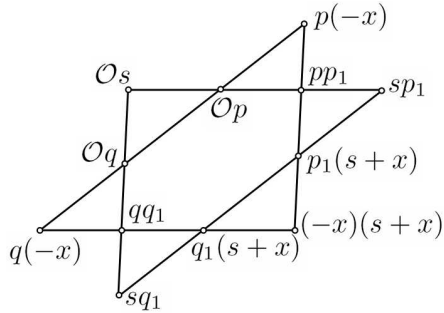}
\caption{The configuration of planes}\label{fig:zvezda}
\end{minipage}
\end{figure}
Here, each plane is denoted by two lines that reflect to each other
on.
\end{proof}

It would be interesting to describe the variety of such
configurations as a moduli-space.

\subsection{General Case}\label{sec:general}

The genus two case we have just considered in detail gives us the
necessary experience of using the billiard constructions to build a
group structure. Moreover, it provides us with the case $n=2$ in the
general construction we are going to present now.

\subsubsection*{Billiard Trajectories and Effective Divisors}

Proposition \ref{prop:elln} is going to serve as the motivation for
defining the operation in the set $\mathcal A_{\ell}$ for
higher-dimensional. Thus, let us now describe in detail the
construction that preceded this proposition in Section
\ref{sec:genus2}.

\smallskip

Suppose that quadrics $\mathcal Q_1$, \dots, $\mathcal Q_n$ from the
confocal family in $\mathbf E^d$ are given. Let $\mathcal O=\ell_0$,
$\ell_1$, \dots, $\ell_n$ be lines in $\mathcal A_{\ell}$ such that
each pair of successive lines $\ell_i$, $\ell_{i+1}$ satisfies the
billiard reflection law at the quadric $\mathcal Q_{i+1}$ ($0\le
i\le n-1$); thus the lines form billiard trajectory $\mathbf
t=(\ell_0,\dots, \ell_n)$.

\smallskip

Let us, for $n\ge1$, define lines
 $\ell_n^{(n)}$, $\ell_n^{(n-1)}$, \dots, $\ell_n^{(1)}$
by the procedure, as follows:

\smallskip

1. Set $\ell_n^{(n)}=\ell_n$.

\smallskip

2. $\ell_n^{(k)}$, for $n-1\ge k\ge1$ is the unique line that
constitutes a DRC with $\ell_{k-1}$, $\ell_{k}$ and
$\ell_n^{(k+1)}$.

\smallskip

Let us note that, in this way, each line in the sequence
$$
\ell_0,\ \ell_1,\ \dots,\ \ell_k,\ \ell_n^{(k+1)},\ \dots,\
\ell_n^{(n)}
$$
is obtained from the previous one by the reflection from
$$
\mathcal Q_1,\ \dots,\ \mathcal Q_k,\ \mathcal Q_n,\ \mathcal
 Q_{k+1},\ \dots,\ \mathcal Q_{n-1}
$$
respectively.

\smallskip

Considering only the initial subsequences $\ell_0$, \dots, $\ell_k$
$(1\le k\le n)$, we may in the same way define lines $\ell_k^{(1)}$,
\dots, $\ell_k^{(k)}$. Notice that, for each $k$, $\ell_k^{(1)}$
intersects $\ell_0$, and these two lines obey the reflection law on
$\mathcal Q_k$.

\smallskip

Thus, we have constructed a mapping $\mathcal D$ from the set
$\mathcal {TB}(\mathcal O)$ of billiard trajectories with the fixed
initial line $\ell_0=\mathcal O$ to the ordered sets of lines from
$\mathcal A_{\ell}$ which intersect $\mathcal O$:
$$
\mathcal D\ :\ \mathbf t=(\ell_0,\dots, \ell_n)\ \mapsto\
(\ell_1^{(1)},\dots, \ell_n^{(1)}),
$$
where line $\ell_k^{(1)}$ intersects $\mathcal O$ according to the
billiard law on the quadric $\mathcal Q_k$.

\smallskip

Doing the opposite procedure, we define morphism $\mathcal B$, the
inverse of $\mathcal D$, which assigns to an $n$-tuple of lines
intersecting $\mathcal O$ the unique billiard trajectory of length
$n$ with $\ell_0=\mathcal O$ as the initial line:
$$
\mathcal B\ :\ (\ell_1^{(1)},\dots, \ell_n^{(1)})\ \mapsto\
(\ell_0,\dots, \ell_n).
$$
Hence, the mapping $\mathcal B$ gives {\em the billiard
representation} of an ordered set of lines intersecting $\mathcal
O$.

\smallskip

In order to consider just divisors of the curve $\mathcal
C_{\mathcal O}$, instead of ordered $n$-tuples of lines, we need to
introduce the following relation $\alpha$ between billiard
trajectories: we say that two billiard trajectories are {\em
$\alpha$-equivalent} if one can be obtained from the other by a
finite set of {\em double reflection moves}. The double reflection
move transforms a trajectory
 $p_1p_2\dots p_{k-1} p_k p_{k+1}\dots p_n$ into trajectory
 $p_1p_2\dots p_{k-1} p_k' p_{k+1}\dots p_n$
if the lines $p_{k-1}, p_k, p_k', p_{k+1}$ form a double reflection
configuration.

\smallskip

It follows directly from The Double Reflection Theorem that two
trajectories $\mathbf t_1,\mathbf t_2$ with the same initial segment
$\mathcal O=\ell_0$ are $\alpha$-equivalent if and only the
$n$-tuples $\mathcal D(\mathbf t_1)$ and $\mathcal D(\mathbf t_2)$
may be obtained from each other by a permutation.

\smallskip

Thus, the mapping $\mathcal D$ may be considered as a mapping from
$$
\widehat{\mathcal{TB}}(\mathcal O)=\mathcal{TB}(\mathcal O)/\alpha
$$
to the set of positive divisors on $\mathcal C_{\mathcal O}$ and it
represents {\em the divisor representation} of billiard
trajectories.

\smallskip

The following lemma, which follows from \cite{Don}, is a sort of
Riemann-Roch theorem in this approach.

\begin{lemma}\label{lema:dijagonala}
A minimal billiard trajectory of length $s$ from $x$ to $y$ is
unique, up the relation $\alpha$. If there are two non
$\alpha$-equivalent trajectories of the same length $k>s$ from $x$
to $y$, then there is a trajectory from $x$ to $y$ of length $k-2$.
\end{lemma}

Now, we are able to introduce an operation, {\it summation}, in the
set $\widehat{\mathcal {TB}}(\mathcal O)$. Given two billiard
trajectories $\mathbf t_1, \mathbf t_2\in
\widehat{\mathcal{TB}}(\mathcal O)$, we define their sum by the
equation:
$$
\mathbf t_1\oplus\mathbf t_2:=\mathcal B(\mathcal D(\mathbf
t_1)+\mathcal D(\mathbf t_2)).
$$
This operation is associative and commutative, according to The
Double Reflection Theorem.

\begin{theorem}
The set $\widehat{\mathcal {TB}}(\mathcal O)=\mathcal {TB}(\mathcal
O)/\alpha$ with the summation operation is a commutative semigroup.
This semigroup is isomorphic to the semigroup of all effective
divisors on the curve $\mathcal C_{\mathcal O}$ that do not contain
the points corresponding to the caustics.
\end{theorem}

Note that in $\widehat{\mathcal {TB}}(\mathcal O)$, the trajectory
consisting only of the single line $\mathcal O$ is a neutral
element.

\smallskip

Define now the following equivalence relation in the set of all
finite billiard trajectories.

\begin{definition}
Two billiard trajectories are {\em $\beta$-equivalent} if they have
common initial and final segments, and they are of the same length.
\end{definition}

The set of all classes of $\beta$-equivalent billiard trajectories
of length $n$, with the fixed initial segment $\mathcal O$, we
denote by $\widetilde{\mathcal{TB}}(\mathcal O)(n)$, and
$$
\widetilde{\mathcal{TB}}(\mathcal O)
=\bigcup_n\widetilde{\mathcal{TB}}(\mathcal O)(n).
$$

\begin{proposition}\label{prop:beta+}
The relation $\beta$ is compatible with the addition of billiard
trajectories.
\end{proposition}

To prove this, we will need the following important lemma:

\begin{lemma}\label{lema:par.nepar}
Let $\mathbf t=(\ell_1,\ell_2,\dots,\ell_{2k},\ell_{2k+1}=\ell_1)$
be a closed billiard trajectory, and $p_1\in\mathcal A_{\ell}$ a
line that intersects $\ell_1$. Construct iteratively lines
$p_2,\dots,p_{2k+1}$ such that quadruples
$p_i,\ell_i,p_{i+1},\ell_{i+1}$ $(1\le i\le 2k)$ form double
reflection configurations. Then $p_{2k+1}=p_1$.
\end{lemma}

\begin{proof}
We are going to proceed with the induction. For $k=2$, the lines
lines $\ell_1,\ell_2,\ell_3,\ell_4$ form a double reflection
configuration, and the statement follows by The Double Reflection
Theorem.

\smallskip

Now, suppose $k>2$. $(\ell_1,\dots,\ell_k)$ and
$(\ell_1=\ell_{2k+1},\ell_{2k},\dots,\ell_k)$ are two billiard
trajectories of the same length from $\ell_1$ to $\ell_k$. If they
are $\alpha$-equivalent, then the claim follows by The Double
Reflection Theorem. If they are not $\alpha$-equivalent, then, by
Lemma \ref{lema:dijagonala}, there is a billiard trajectory
$\ell_1'=\ell_1,\ell_2',\dots,\ell_{k-2}'=\ell_k$. The statement now
follows from the inductive hypothesis applied to trajectories
$$(\ell_1,\ell_2,\dots,\ell_k=\ell_{k-2}',\ell_{k-3}',\dots,\ell_1'=\ell_1)$$
and
$$(\ell_1',\ell_2',\dots,\ell_{k-2}'=\ell_k,\ell_{k+1},\dots,\ell_{2k+1}=\ell_1).$$
\end{proof}

\noindent{\it Proof of Proposition \ref{prop:beta+}}. We need to
prove the following:
$$
\mathbf t_1\sim_{\beta}\mathbf t_1',\ \mathbf t_2\sim_{\beta}\mathbf
t_2'\ \Rightarrow\ \mathbf t_1+\mathbf t_2\sim_{\beta}\mathbf
t_1'+\mathbf t_2'.
$$
Clearly, it is enough to prove this relation for the case when
$\mathbf t_2=\mathbf t_2'$ and the length of $\mathbf t_2$ is equal
to $2$. Suppose that $\mathbf t_2=(\mathcal O, p)$, where $p$ is
obtained from $\mathcal O$ by the reflection on quadric $\mathcal
Q_p$. Then trajectories $\mathbf t_1+\mathbf t_2$, $\mathbf
t_1+\mathbf t_2'$ are obtained by adding one segment to $\mathbf
t_1$, $\mathbf t_1'$ respectively. These segments satisfy the
reflection law on $\mathcal Q_p$ with the final segments of $\mathbf
t_1$, $\mathbf t_1'$. Since $\mathbf t_1\sim_{\beta}\mathbf t_1'$,
their final segments coincide and the statement follows from Lemma
\ref{lema:par.nepar} applied to the trajectory
$$
 (\ell_1, \ell_2, \dots, \ell_n=\ell_n', \ell_{n-1}', \ell_{n-2}', \dots, \ell_1'=\ell_1),
$$
where $\mathbf t_1=(\ell_1,\dots,\ell_n)$, $\mathbf
t_2=(\ell_1',\dots,\ell_n')$.
 \hfill\hfill$\Box$\linebreak

\subsubsection*{The Group Structure in $\mathcal A_{\ell}$}

We wish to use the constructed algebra on the set of billiard
trajectories, in order to obtain an algebraic structure on $\mathcal
A_{\ell}$, such that the operation is naturally connected with the
billiard reflection law.

\smallskip

From \cite{Tyu} we get

\begin{theorem}\label{teorema:dostiznost}
For any two given lines $x$ and $y$ from $\mathcal A_{\ell}$, there
is a system of at most $d-1$ quadrics from the confocal family, such
that the line $y$ is obtained from $x$ by consecutive reflections on
these quadrics.
\end{theorem}

The divisor representation of the corresponding billiard trajectory
of length $s\le d-1$ will be called {\em the $s$-brush of $y$
related to $x$}.

\smallskip

Now, we can define a group structure in $\mathcal A_{\ell}$
associated with a fixed line in this set.

\medskip

\noindent{\em Neutral element}

\smallskip

Let us fix a line $\mathcal O\in\mathcal A_{\ell}.$

\medskip

\noindent{\em Inverse element}

\smallskip

Let $x$ be an arbitrary line in $\mathcal A_{\ell}$, and $\mathcal
D(x)$ the divisor representation of the minimal billiard trajectory
connecting $\mathcal O$ with $x$. Define $-x$ as the final segment
of the billiard trajectory $\mathcal B(\tau\mathcal D(x))$, where
$\tau$ is the hyperelliptic involution of $\mathcal C_{\mathcal O}$.

\medskip

\noindent{\em Addition}

\smallskip

For two lines $x$ and $y$ from $\mathcal A_{\ell}$, denote their
brushes related to $\mathcal O$ as $S_1$ and $S_2$. Define
$$
x+y:=
(-1)^{|S_1|+|S_2|+1}\mathcal{EB}(\tau^{|S_1|+1}(S_1),\tau^{|S_2|+1}(S_2)),
$$
where $\mathcal{EB}(\tau^{|S_1|+1}(S_1),\tau^{|S_2|+1}(S_2))$ is the
final segment of the billiard trajectory
$\mathcal{B}(\tau^{|S_1|+1}(S_1),\tau^{|S_2|+1}(S_2))$.

\smallskip

From all above, similarly as in the genus $2$ case, we get

\begin{theorem}
The set $\mathcal A_{\ell}$ with the operation defined above is an
Abelian group.
\end{theorem}

Let us observe that another kind of divisor representation of
billiard trajectories may be introduced, that associates a divisor
of degree $0$ to a given trajectory. Such a representation was used
in \cite{Dar3} and explicitly described in \cite{DR3, DR4}. The
positive part of this representation coincide with the divisor
obtained by $\mathcal B$, and the negative part is invariant for the
hyperelliptic involution. Moreover, as it follows from \cite{Dar3},
two trajectories with the initial segment $\mathcal O$ will have the
same final segments if and only if their representations are
equivalent divisors. Thus, it follows that the group structure on
$\mathcal A_{\ell}$ is isomorphic to the quotient of the Jacobian of
the curve $\mathcal C_{\mathcal O}$ by a finite subgroup which is
generated by the points corresponding to the caustic quadrics.

\section{Billiard Algebra, Theorems of Poncelet Type and Their Generalizations}
\label{sec:algebra}

\subsection{Billiard Algebra, Weak Poncelet Trajectories and
Theorems of Poncelet-Cayley's Type}\label{sec:weak}

Now, as a consequence of Theorem \ref{teorema:dostiznost} we are
able to introduce the following hierarchies of notions.

\begin{definition}\label{definicija:s-mimoilaznost}
For  two given lines $x$ and $y$ from $\mathcal A_{\ell}$ we say
that they are {\em $s$-skew} if $s$ is the smallest number such that
there exist a system of $s+1\le d-1$ quadrics $\mathcal Q_k$,
$k=1,...,s+1$ from the confocal family, such that the line $y$ is
obtained from $x$ by consecutive reflections on $\mathcal Q_k$. If
the lines $x$ and $y$ intersect, they are {\em $0$-skew}. They are
{\em $(-1)$-skew} if they coincide.
\end{definition}

\begin{definition}\label{definicija:s-Poncelet}
Suppose that a system $S$ of $n$ quadrics $\mathcal Q_1$, \dots,
$\mathcal Q_n$ from the confocal family is given. For a system of
lines $\mathcal O_0$, $\mathcal O_1$, \dots, $\mathcal O_n$  in
$\mathcal A_{\ell}$ such that each pair of successive lines
$\mathcal O_i$, $\mathcal O_{i+1}$ satisfies the billiard reflection
law at $\mathcal Q_{i+1}$ $(0\le i\le n-1)$, we say that it forms an
{\em $s$-weak Poncelet trajectory of length $n$ associated to the
system $S$} if the lines $\mathcal O_0$ and $\mathcal O_n$ are
$s$-skew.
\end{definition}

For $s$-weak Poncelet trajectories we will also sometimes say {\em
$(d-s-2)$-resonant billiard trajectories}. Periodic trajectories or
generalized classical Poncelet polygons are $(-1)$-weak Poncelet
trajectories or, in other words, $(d-1)$-resonant billiard
trajectories, and they are described analytically in \cite{DR4},
\cite{DR3}.

\smallskip

Our next goal is to get complete analytical description of $s$-weak
Poncelet trajectories of length $r$, generalizing in such a way the
results from \cite{DR4}. Here, we are going to use fully the tools
and the power of billiard algebra.

\smallskip

To fix the idea, let us consider first the system $S$ consisting on
$r$ equal quadrics $\mathcal Q_1= \dots = \mathcal Q_r$ from the
confocal family. Suppose a system of lines $\mathcal O_0$, $\mathcal
O_1$, \dots, $\mathcal O_r$  in $\mathcal A_{\ell}$ forms an
$s$-weak Poncelet trajectory of length $r$ associated to the system
$S$. Then
$$
\mathcal O_r = r\mathcal O_1^{(1)},
$$
with some line $\mathcal O_1^{(1)}$ which intersects $\mathcal O_0$.
Again, from the condition that $\mathcal O_r$ and $\mathcal O_0$ are
$s$-skew we get
$$\mathcal O_r =
\mathcal O_1^{'(1)}+\dots+\mathcal O_{s+1}^{'(1)}$$ with some
lines $P_i=\mathcal O_i^{'(1)}$ which intersect $\mathcal O_0$.
From the last two equations we come to the conclusion

\begin{proposition}\label{propozicija:kejli1}
The existence of an $s$-weak Poncelet trajectory of length $r$ is
equivalent to existence of a meromorphic function $f$ on the
hyperelliptic curve $\mathcal C_{\ell}$ such that $f$ has a zero of
order $r$ at $P=\mathcal O_1^{(1)}$, a unique pole at ``infinity"
$E$, and the order of the pole is equal to $n=r+s+1$.
\end{proposition}

\smallskip

Now, we are going to derive explicit analytical condition of
Cayley's type. As in \cite{DR4} we consider the space $\mathcal
L(nE)$ of all meromorphic functions on $\mathcal C_{\ell}$ with the
unique pole at the infinity point $E$ of the order not exceeding
$n$. Let $(f_1,\dots, f_k)$ be one of the bases of this space,
$k=\dim \mathcal L(nE)$. Consider the vectors
$$v_1,\dots,v_r\in \mathbf C^{k},$$
where $v_i^j=f_j^{(i-1)}(P)$ and vectors
$$u_1,\dots,u_{s+1}\in
\mathbf C^{k},$$
 with $u_i^j=f_j(P_i)$.
From the condition (see \cite {DR4})
$$\rank[v_1,\dots,v_r,u_1,\dots,u_{s+1}]<n-g+1$$
we get the condition
$$\rank[v_1,\dots,v_r]<r+s-g+2=r+s-d+3.$$

\smallskip

Now, we can rewrite it in the  form usual for the Cayley's type
conditions.

\begin{theorem}\label{teorema:kejlijevi.uslovi}
The existence of an $s$-weak Poncelet trajectory of length $r$ is
equivalent to:
$$
 \rank\left(\begin{array}{llll}
B_{d+1} & B_{d+2} & \dots & B_{m+1}\\
B_{d+2} & B_{d+3} & \dots & B_{m+2}\\
\dots & \dots & \dots & \dots\\
B_{d+m-s-2} & B_{d+m-s-1} & \dots & B_{r-1}
\end{array}\right)<m-d+1,
$$
when $r+s+1=2m$, and
$$
 \rank\left(\begin{array}{llll}
B_{d} & B_{d+1} & \dots & B_{m+1}\\
B_{d+1} & B_{d+2} & \dots & B_{m+2}\\
\dots & \dots & \dots & \dots\\
B_{d+m-s-2} & B_{d+m-s-1} & \dots & B_{r-1}
\end{array}\right)<m-d+2,
$$
when $r+s+1=2m+1$.

With $B_0,B_1,B_2,\dots$, we denoted the coefficients in the Taylor
expansion of function $y=\sqrt{\mathcal P(x)}$ in a neighbourhood of
$P$, where $y^2=\mathcal P(x)$ is the equation of the generalized
Cayley curve, with the polynomial $\mathcal P(x)$ given by
(\ref{eq:polinom}).
\end{theorem}

\begin{proof}
Denote by $\mathcal L((r+s+1)E)$ the linear space of all meromorphic
functions on $\mathcal C_{\ell}$, having a unique pole at the
infinity point $E$, with the order not exceeding $r+s+1$. By the
Riemann-Roch theorem:
\begin{itemize}
\item[$(i)$]
$\dim\mathcal L((r+s+1)E)=[\frac{r+s+1}2]+1$ if $r+s+1\le2d-1$,

\item[$(ii)$]
$\dim\mathcal L((r+s+1)E)=r+s-g+1$ if $r+s+1$ is even and greater
than $2d-2$,

\item[$(iii)$]
$\dim\mathcal L((r+s+1)E)=r+s-g+2$ if $r+s+1$ is even and greater
than $2d-2$.
\end{itemize}
We may choose the following bases in each of the three cases:
\begin{itemize}
\item[$(i)$]
$1,x,\dots,x^m$, where $m=[\frac{r+s+1}2]\le d$;

\item[$(ii)$]
$1,x,\dots,x^m, y, xy, \dots, x^{m-d}$, for $r+s+1=2m\ge 2d-2$;

\item[$(iii)$]
$1,x,\dots,x^m, y, xy, \dots, x^{m-d+1}$, for $r+s+1=2m+1\ge 2d-2$.
\end{itemize}
Now, the statement follows from the considerations preceding the
theorem.
\end{proof}

\smallskip

\begin{example}
For $s=-1$, the inequalities in Theorem
\ref{teorema:kejlijevi.uslovi} become the conditions for periodic
billiard trajectories (see \cite{DR3,DR4}).
\end{example}

\begin{example}
The condition for existence of a $(d-3)$-weak Poncelet trajectory of
length $r$ is equivalent to:
$$
 \det\left(\begin{array}{llll}
B_{d+1} & B_{d+2} & \dots & B_{m+1}\\
B_{d+2} & B_{d+3} & \dots & B_{m+2}\\
\dots & \dots & \dots & \dots\\
B_{m+1} & B_{m+2} & \dots & B_{r-1}
\end{array}\right)=0,\quad r+d-2=2m;
$$
$$
 \det\left(\begin{array}{llll}
B_{d} & B_{d+1} & \dots & B_{m+1}\\
B_{d+1} & B_{d+2} & \dots & B_{m+2}\\
\dots & \dots & \dots & \dots\\
B_{m+1} & B_{m+2} & \dots & B_{r-1}
\end{array}\right)=0,\quad r+d-2=2m+1.
$$
\end{example}

\begin{example}
If $r+s<2d-1$, then, as it follows from the proof of Theorem
\ref{teorema:kejlijevi.uslovi}, an $s$-weak Poncelet trajectory of
length $r$ may exist only if the hyper-elliptic curve $\mathcal
C_{\ell}$ is singular.
\end{example}

\subsection{Remark on Generalized Weyr's Theorem and Griffiths-Harris Space
Poncelet Theorem in Higher Dimensions}\label{sekcija:wgh}

In this section we obtain higher dimensional generalizations of the
results of \cite {We}, \cite{Hu}, \cite{GH1}. Nice exposition of
those classical results one can find in \cite{BB}. The dual version
of \cite{GH1} is given in \cite{CCS}.

\smallskip

Each quadric $\mathcal Q$ in $\mathbf P^{2d-1}$ contains at most two
unirational families of $(d-1)$-dimensional linear subspaces. Such
unirational families are usually called {\em rulings of the
quadric}.

\smallskip

\begin{theorem}\label{th:gen.weyr}
Let $\mathcal Q_1$, $\mathcal Q_2$ be two general quadrics in
$\mathbf P^{2d-1}$ with the smooth intersection $V$ and $\mathcal
R_1$, $\mathcal R_2$ their rulings. If there exists a closed chain
$$
L_1,\ L_2,\ \dots,\ L_{2n},\ L_{2n+1}=L_1
$$
of distinct $(d-1)$-di\-men\-sional linear subspaces, such that
$L_{2i-1}\in\mathcal R_1$, $L_{2i}\in\mathcal R_2$ $(1\le i\le n)$
and $L_j\cap L_{j+1}\in F(V)$ $(1\le j\le 2n)$, then there are such
closed chains of subspaces of length $2n$ through any point of
$F(V)$.
\end{theorem}

\begin{proof}
Each of the unirational families $\mathcal R_i$ determines an
involution $\tau_i$ on Abelian variety $F(V)$. Such an involution
interchanges two $(d-2)$-intersections of an element of $\mathcal
R_i$ with $V$. Denote by $\tr:F(V)\to F(V)$ their composition and by
$L:=L_{2n}\cap L_{1}\in F(V)$. Since $\tr$ is a translation on
$F(V)$ satisfying $\tr^n(L)=L$ we see that $\tr$ is of order $n$ and
the theorem follows.
\end{proof}

\begin{definition}
We will call the chains considered in Theorem \ref{th:gen.weyr} {\em
generalized Weyr's chains}.
\end{definition}

The theorem can be adjusted for nonsmooth intersections, but we are
not going into details. Instead, we consider the case of two
quadrics (\ref{eq:kvadrika1}) and (\ref{eq:kvadrika2}) as in Section
\ref{sekcija:morfizam}.  By using the projection $\pi^*$ we get:

\begin{proposition}\label{propozicija:weyr.poncelet}
A generalized Weyr chain of length $2n$ projects into a Poncelet
polygon of length $2n$ circumscribing the quadrics $\mathcal
Q^p_{\alpha_1},\dots,\mathcal Q^p_{\alpha_{d-1}}$ and alternately
inscribed into two fixed confocal quadrics (projections of $\mathcal
Q_1, \mathcal Q_2$). Conversely, any such a Poncelet polygon of the
length $2n$ circumscribing the quadrics $\mathcal
Q^p_{\alpha_1},\dots,\mathcal Q^p_{\alpha_{d-1}}$ and alternately
inscribed into two fixed confocal quadrics can be lifted to a
generalized Weyr chain of length $2n$.
\end{proposition}

\begin{proof}
It follows from Lemma 4.1 and Corollary 4.2 of \cite{Kn} and Lemma
\ref{lema:veza}.
\end{proof}

Thus, we obtained in  a correspondence between generalized Weyr
chains and Poncelet polygons subscribed in $d-1$ given quadrics and
alternately inscribed in two quadrics from some confocal family.
Such Poncelet polygons have been completely analytically described,
among others, in \cite{DR4} (see Example 4 from there).

\smallskip

Let us note that correspondence between classical Weyr's theorem in
$\mathbf P^3$ and classical Poncelet theorem about two conics in a
plane was observed by Hurwitz in \cite{Hu}. Nevertheless the
projection we used here is not a straightforward generalization of
the one used by Hurwitz.

\smallskip

Polygonal lines, circumscribed about one conic and alternately
inscribed in two conics, appear as an example in \cite{Ves}. Our
result, applied to the lowest dimension, gives the condition for the
closeness of such a polygonal line, see Corollary 1 from \cite{DR4}.

\smallskip

Let us also mention that if rulings $\mathcal R_1$ and $\mathcal
R_2$ are connected by a generalized Weyr's chain of length $2n$,
then the same is true for $\mathcal R_2$, $\mathcal R_1$ and also
for the pair $\mathcal R'_1$, $\mathcal R'_2$ of the complementary
rulings of $\mathcal Q_1$ and $\mathcal Q_2$.

\smallskip

Now we are able to present a new higher-dimensional generalization
of the Griffiths-Harris Space Poncelet theorem from \cite{GH1}.

\smallskip

\begin{theorem}\label{teorema:uopstena.grifits.haris}
Let $\mathcal Q^*_1$ and $\mathcal Q^*_2$ be the duals of two
general quadrics in $\mathbf P^{2d-1}$ with the smooth intersection
$V$. Denote by $\mathcal R_i,\mathcal R'_i$ pairs of unirational
families of $(d-1)$-dimensional subspaces of $\mathcal Q^*_i$.
Suppose there are generalized Weyr's chains between $\mathcal R_1$
and $\mathcal R_2$ and between $\mathcal R_1$ and $\mathcal R'_2$.
Then there is a finite polyhedron inscribed and subscribed in both
quadrics $\mathcal Q_1$ and $\mathcal Q_2$. There are infinitely
many such polyhedra.
\end{theorem}

\smallskip

The polyhedra from the previous theorem can be described in more
details as arrangements of $d$-dimensional ``faces" which are
bitangents of the quadrics $\mathcal Q_1$ and $\mathcal Q_2$. Their
intersections are $(d-1)$-dimensional ``edges" which alternately
belong to the rulings of $\mathcal Q^*_1$ and $\mathcal Q^*_2$. The
intersections of ``edges" are $(d-2)$-dimensional ``vertices" which
belong to $F(V)$.

\smallskip

The proof of the last theorem is based on the following Lemma from
\cite{Don}.

\begin{lemma}\label{lema:odredjenost}
Let $\mathcal Q\subset\mathbf P^{2d-1}$ be a quadric of the rank not
less than $2d-1$ and $x\subset\mathcal Q$ a linear subspace of the
dimension $d-2$. If $\mathcal Q$ is singular, suppose additionally
that $x$ does not contain the vertex. Then, for each ruling
$\mathcal R$ of $\mathcal Q$, there is a unique $(d-1)$-dimensional
linear subspace $s=s(\mathcal R,x)$ such that $s\in\mathcal R$ and
$x\subset s$.
\end{lemma}

\subsection{Poncelet-Darboux Grid and Higher Dimensional
Generalizations}\label{sec:grid}

\subsubsection*{Historical Remark: on Darboux's Heritage}

In \cite{Dar3} (Volume 3, Book VI, Chapter I), Darboux proved that
Liouville surfaces are exactly these having an orthogonal system of
curves that can be regarded in two or, equivalently, infinitely many
different ways, as geodesic conics. These coordinate curves are
analogues of systems of confocal conics in the Euclidian plane.
Knowing this, essential properties of conics may be generalized to
all Liouville surfaces (see \cite{Dar3}, and also \cite{DR4} for a
review and clarifications). Here is the citation of one of these
properties:

\smallskip

{\em ``Considerons un polygone variable dont tous les c\^ot\'es sont
des lignes g\'eo\-d\'e\-si\-ques tangentes a une m\^eme courbe
coordonn\'ee. Le dernier sommet de ce polygone d\'ecrira aussi une
des courbes coordonn\'ees et il en sera de m\^eme de points
d'intersection de deux c\^ot\'es quelconques de ce polygone.''}

\smallskip

This statement is not only the generalization of Poncelet theorem to
Liouville surfaces. Additionally, for a Poncelet polygon
circumscribed about a fixed coordinate curve, with each vertex
moving along one of the coordinate curves, Darboux stated here that
the intersection point of two arbitrary sides of the polygon is also
going to describe a coordinate curve.

\smallskip

Let us note that a weaker version of this claim, although with some
improvements, has been recently rediscovered in \cite{Svarc}, and a
more elementary proof of the result from \cite{Svarc} has been
published in \cite{LeviTab}. The main theorem of \cite{Svarc}
corresponds only to a fixed Poncelet polygon associated to a pair of
ellipses in the Euclidean plane. In a sense, this gives a new
argument in favour of our observation in \cite{DR3} and \cite{DR4},
that the work of Darboux connected with billiards and Poncelet
theorem seems to be unknown to nowadays mathematicians.

\smallskip

This section is devoted to a multi-dimensional generalization of the
Darboux theorem, related to billiard trajectories within an
ellipsoid in the $d$-dimensional Euclidean space.

\subsubsection*{Poncelet-Darboux Grid in Euclidean Plane}

Before starting with higher-dimensional generalizations, let us give
some improvements together with a simpler proof of the statement
corresponding to the plane case.

\begin{theorem}\label{th:gen.grida}
Let $\mathcal E$ be an ellipse in $\mathbf E^2$ and
$(a_m)_{m\in\mathbf Z}$, $(b_m)_{m\in\mathbf Z}$ be two sequences of
the segments of billiard trajectories $\mathcal E$, sharing the same
caustic. Then all the points $a_m\cap b_m$ $(m\in\mathbf Z)$ belong
to one conic $\mathcal K$, confocal with $\mathcal E$.

\smallskip

Moreover, under the additional assumption that the caustic is an
ellipse, we have: if both trajectories are winding in the same
direction about the caustic, then $\mathcal K$ is also an ellipse;
if the trajectories are winding in opposite directions, then
$\mathcal K$ is a hyperbola. (See Figure \ref{fig:grida_elipsa}.)
\begin{figure}[h]
\centering
\begin{minipage}[t]{0.44\textwidth}
\centering
\includegraphics[width=5cm,height=4cm]{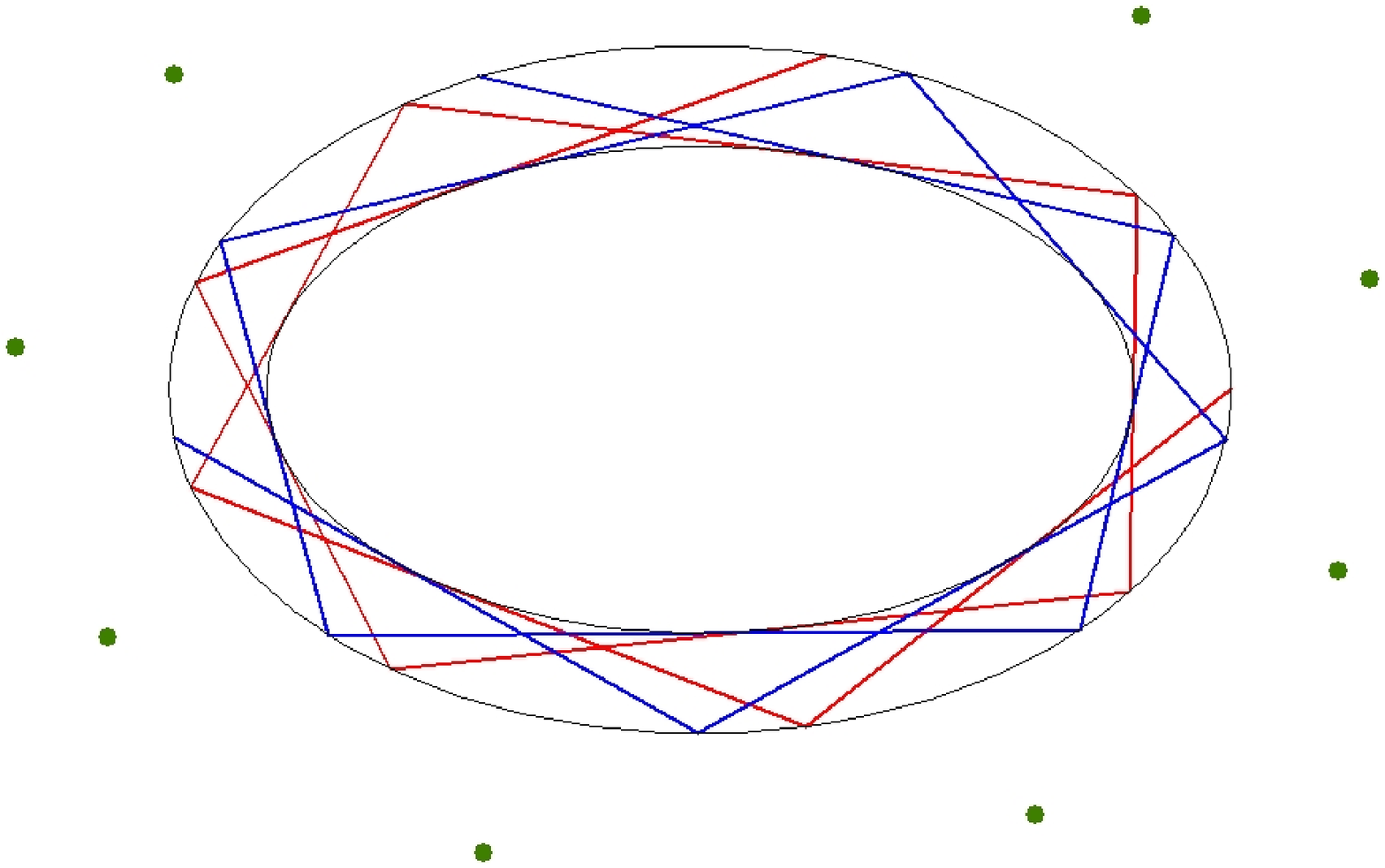}
\end{minipage}
\begin{minipage}[t]{0.44\textwidth}
\centering
\includegraphics[width=5cm,height=4cm]{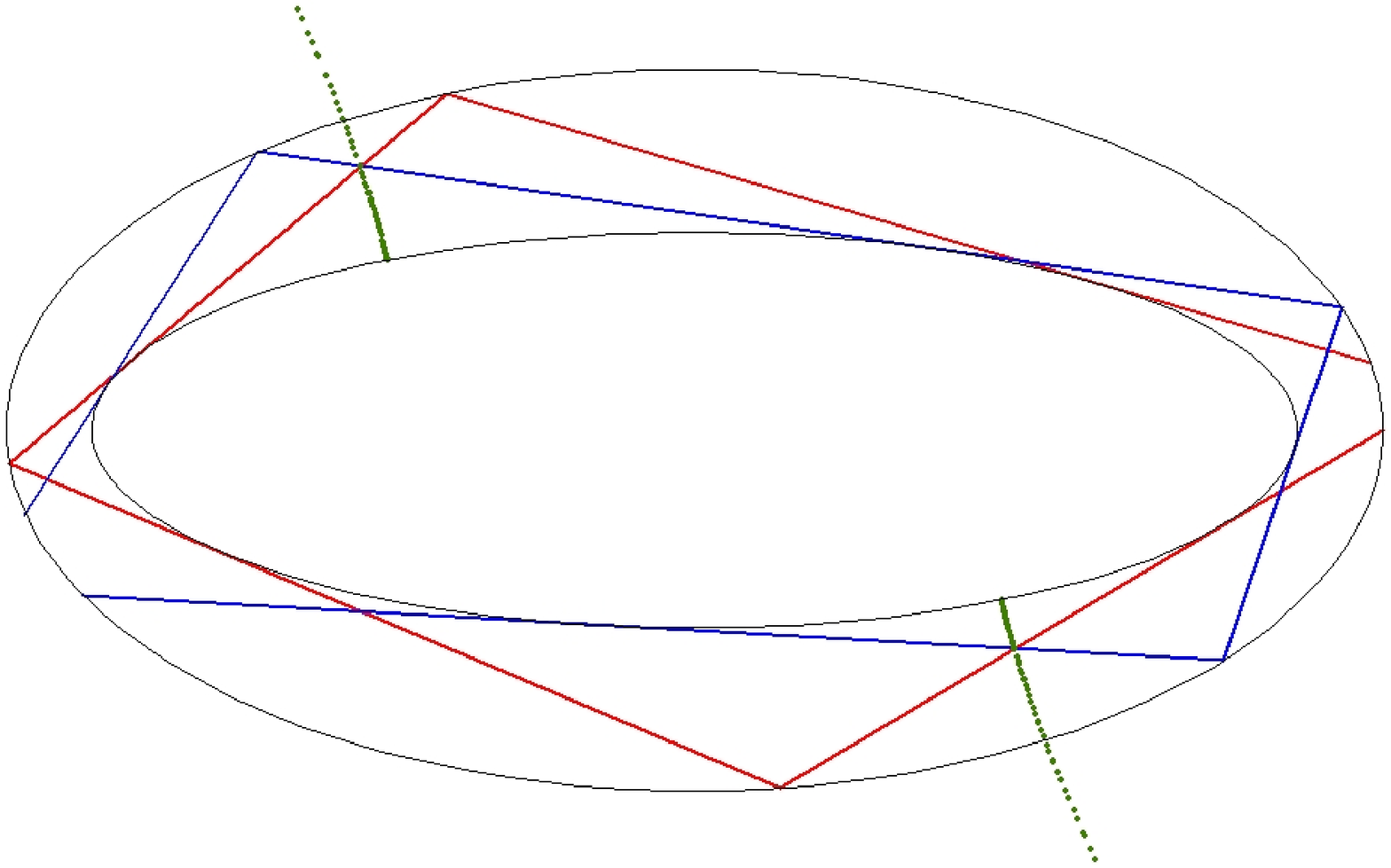}
\end{minipage}
\caption{Billiard trajectories with an ellipse as a caustic and the
intersection points of the corresponding
segments}\label{fig:grida_elipsa}
\end{figure}

\smallskip

For a hyperbola as a caustic, it holds: if segments $a_m$, $b_m$
intersect the long axis of $\mathcal E$ in the same direction, then
$\mathcal K$ is a hyperbola, otherwise it is an ellipse. (See Figure
\ref{fig:grida_hip}.)
\begin{figure}[h]
\centering
\begin{minipage}[t]{0.44\textwidth}
\centering
\includegraphics[width=5cm,height=4cm]{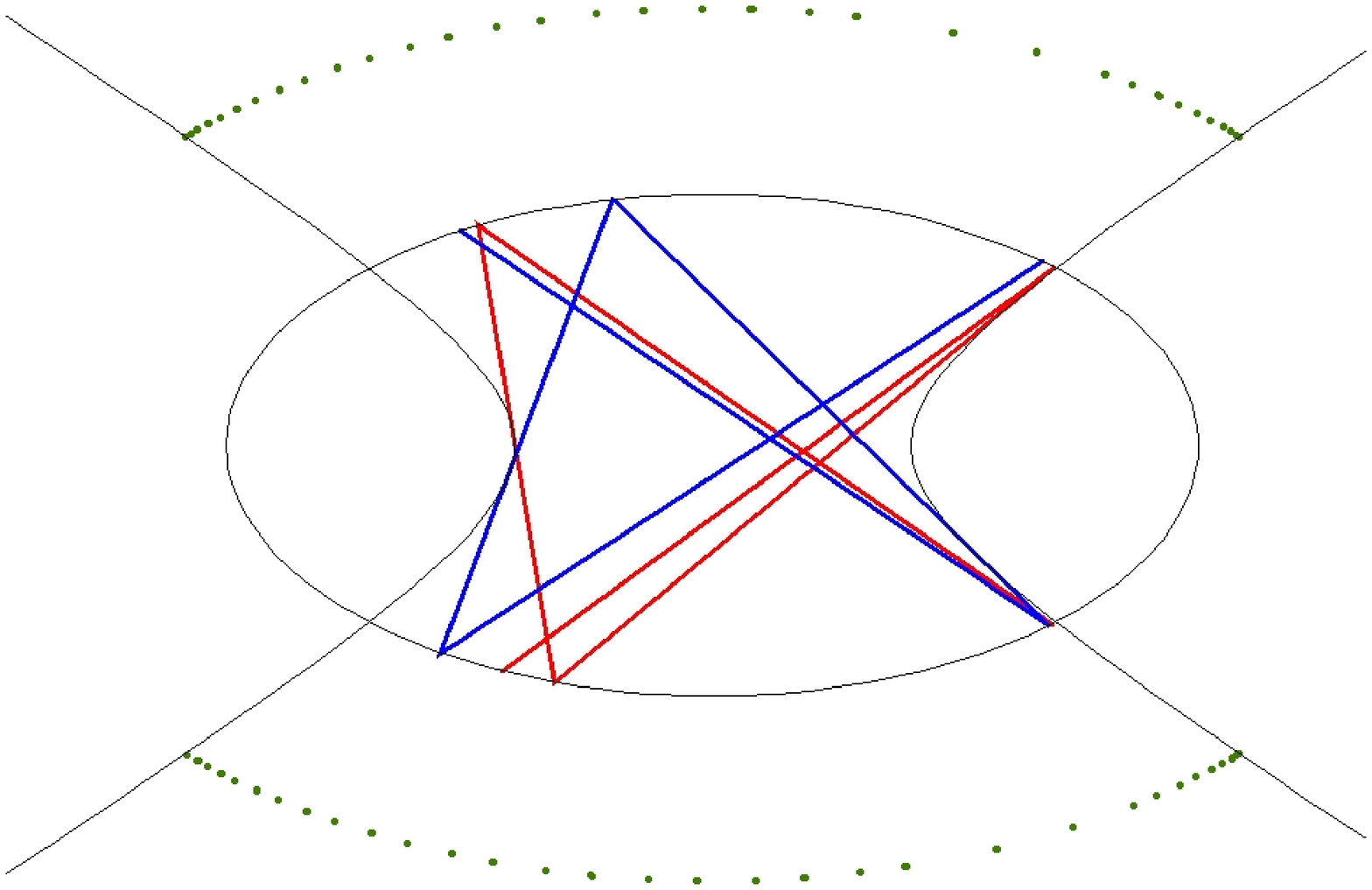}
\end{minipage}
\begin{minipage}[t]{0.44\textwidth}
\centering
\includegraphics[width=5cm,height=4cm]{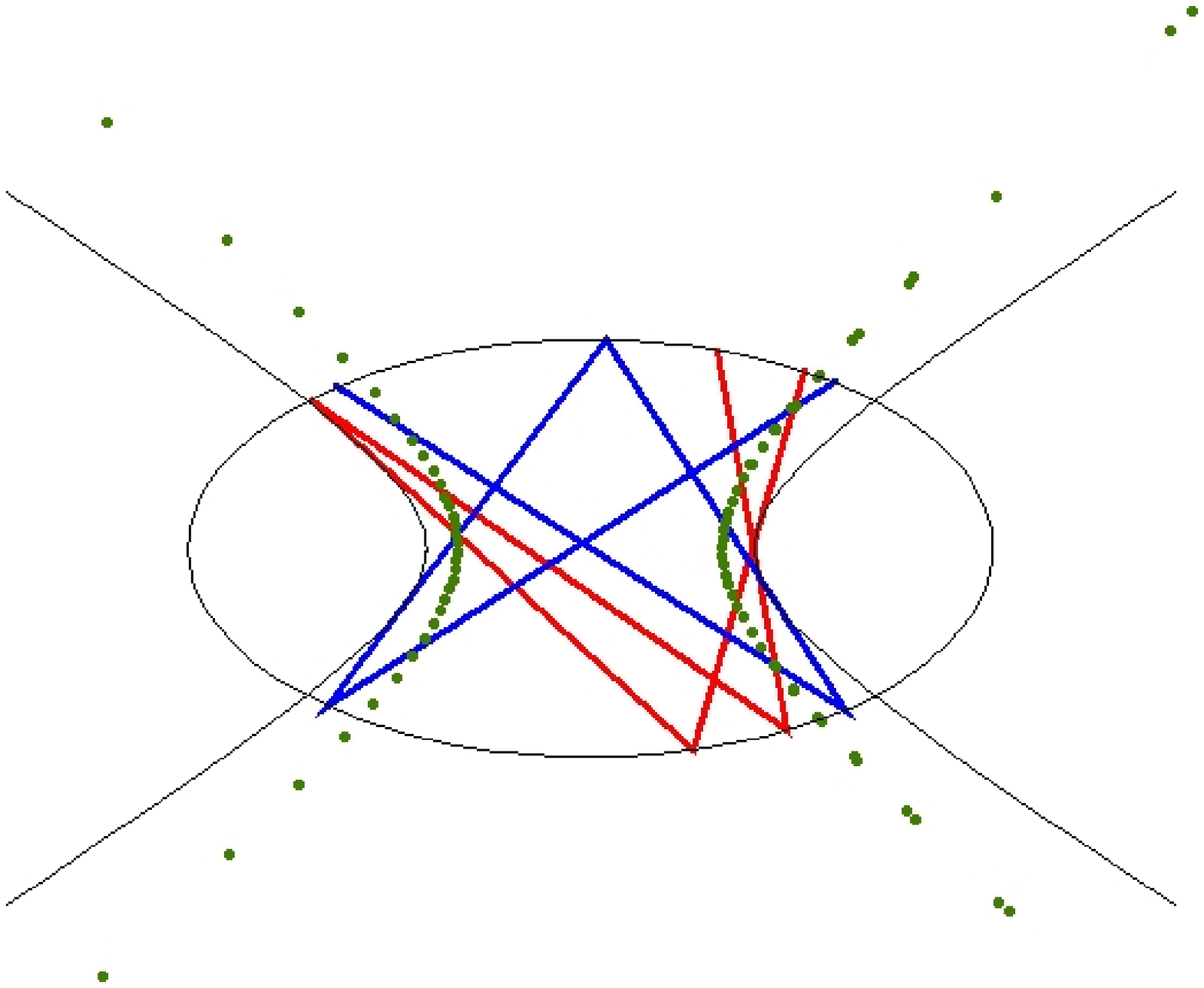}
\end{minipage}
\caption{Billiard trajectories with a hyperbola as a caustic and the
intersection points of the corresponding
segments}\label{fig:grida_hip}
\end{figure}
\end{theorem}

\begin{proof}
The statement follows by application of the Double Reflection
Theorem. Namely, the lines $a_0$ and $b_0$ intersect at a point that
belongs to one ellipse and one hyperbola from the confocal family.
They satisfy the reflection law on exactly one of these two curves,
depending on the orientation of the billiard motion along the lines.
Now, by the Double Reflection Theorem, $a_1$ and $b_1$ satisfy the
reflection law on the same conic. The same is true, for any pair
$a_m$, $b_m$, by the induction.

\smallskip

For the second part of the theorem, it is sufficient to observe that
the winding direction about an ellipse is changed by the reflections
on the hyperbolae, and preserved by the reflections on the ellipses
from the confocal family. If an oriented line is placed between the
foci, then the direction it is intersecting the axis containing the
foci is changed by the reflections of the ellipses and preserved by
the reflections on the hyperbolae.
\end{proof}

\begin{proposition}
Let $(a_m)_{m\in\mathbf Z}$, $(b_m)_{m\in\mathbf
Z}$ be two sequences of the segments of billiard trajectories within
the ellipse $\mathcal E$, sharing the same caustic. If the caustic
is an ellipse and the trajectories are winding in the opposite
directions about it, then all the points $a_m\cap b_m$ $(m\in\mathbf
Z)$ are placed on two centrally symmetric half-branches of a
hyperbola confocal with $\mathcal E$.
\end{proposition}

\begin{proof}
Denote by $\mathcal E_c$ the ellipse which is the common caustic of
the given billiard trajectories. It is possible to introduce a
metric $\mu$ on $\mathcal E_c$, such that $\mu(AB)=\mu(CD)$ if and
only if the tangent lines at $A,B$ and at $C,D$ intersect on the
same ellipse of the confocal family. In particular, this means that
the points where two consecutive billiard segments touch the caustic
$\mathcal E'$ are always at the fixed distance from each other.

\smallskip

Take that $a_0$ and $b_0$ intersect on the upper left half-branch of
the hyperbola $\mathcal K$ (see Figure \ref{fig:grida_elipsa}). Then
it may be proved from The Double Reflection Theorem, that each
touching points of these two segments with caustic $\mathcal E'$ are
at equal $\mu$-distances from the upper left intersection point of
$\mathcal K$ and $\mathcal E_c$. Since $\mu$ is symmetric with
respect to the coordinate centre, they are also equally distanced
from the lower right intersection point, but not from the other two
intersection points. The same holds for any pair $a_m,b_m$. Thus,
their intersections lies on two centrally symmetric half-branches of
the hyperbola.
\end{proof}

\medskip

Now, the generalized claim about Poncelet-Darboux grids is
immediately following from Theorem \ref{th:gen.grida}.

\begin{theorem}\label{th:grida}
Let $(\ell_m)_{m\in\mathbf Z}$ be the sequence of segments of a
billiard trajectory within the ellipse $\mathcal E$. Then each of
the sets
 $\mathrm P_k= \bigcup_{i-j=k}\ell_i\cap\ell_j$, $\mathrm Q_k =
 \bigcup_{i+j=k}\ell_i\cap\ell_j$, $(k\in\mathbf Z)$
belongs to a single conic confocal to $\mathcal E$.

\smallskip

If the caustic of the trajectory $(\ell_m)$ is an ellipse, then the
sets $\mathrm P_k$ are placed on ellipses and $\mathrm Q_k$ on
hyperbolae. If the caustic is a hyperbola, then the sets $\mathrm
P_k$, $\mathrm Q_k$ are placed on ellipses for $k$ even and on
hyperbolae for $k$ odd.
\end{theorem}

\begin{proof}
To show the statement for $\mathrm P_k$, take $a_m=\ell_m$,
$b_m=\ell_{m+k}$ and apply Theorem \ref{th:gen.grida}. For $\mathrm
Q_k$, take $a_m=\ell_m$, $b_m=\ell_{k-m}$.
\end{proof}

\begin{remark}
Notice that Theorem \ref{th:grida} is more general then the one
given in \cite{Svarc, LeviTab}, since we do not suppose that the
billiard trajectory is closed. Also, the statement can be
formulated for an arbitrary conic, not only for an ellipse.
\end{remark}

The main statement proved in \cite{Svarc} and \cite{LeviTab} is a
special case consequence of Theorem \ref{th:grida}:

\begin{corollary}[\cite{Svarc}, \cite{LeviTab}]
Let $(\ell_m)$ be a closed billiard trajectory within an ellipse,
with the elliptical caustic. Each set $\mathrm P_k$ lie on an
ellipse confocal to $\mathcal E$, and $\mathrm Q_k$ on a confocal
hyperbola. {\em (See Figure \ref{fig:grida_elipsa}).}
\end{corollary}

Let us show one interesting property of Poncelet-Darboux grids.

\begin{proposition}
Let $(\ell_m)$ be a billiard trajectory within ellipse $\mathcal E$,
with the elliptical caustic $\mathcal E_c$. Then the ellipse
containing the set $\mathrm P_k$ depends only on $k$, $\mathcal E$
and $\mathcal E_c$. In other words, this ellipse will remain the
same for any choice of billiard trajectory within $\mathcal E$ with
the caustic $\mathcal E_c$.
\end{proposition}

\begin{proof}
This claim may be proved by the use of the Double Reflection
Theorem, similarly as in Theorem \ref{th:gen.grida}. Nevertheless,
we are going to show it in another way, that gives a possibility of
explicit calculation of the ellipse containing $\mathrm P_k$.

\smallskip

To the billiard within $\mathcal E$ with the fixed caustic $\mathcal
E_c$, we may associate an elliptic curve (see \cite{Leb,GH2,Mo,MV}).
Each point of the curve corresponds to a conic from the confocal
family. On the other hand, the billiard motion may be viewed as the
linear motion on the Jacobian of the curve (i.e.\ on the curve
itself, since this is the elliptic case), with translation jumps
corresponding to reflections. More precisely, the translation is
exactly by this value on the elliptic curve, which is associated to
the ellipse that a segment is reflected on.

\smallskip

Thus, if the point $M$ on the elliptic curve is associated to the
ellipse $\mathcal E$, then the set $\mathrm P_k$ is placed on the
ellipse corresponding to $kM$.
\end{proof}

\subsubsection*{Grids in Arbitrary Dimension}

Although Theorem \ref{th:gen.grida} gives certain progress in
understanding of Poncelet-Dar\-boux grids, the essential
breakthrough in this matter represents, in our opinion, the study of
the higher-dimensional situation. This analysis is based on
introduction of higher-dimensional analogues of grids and our notion
of $s$-skew lines.

\begin{theorem}\label{th:gen.grida.d}
Let $(a_m)_{m\in\mathbf Z}$, $(b_m)_{m\in\mathbf Z}$ be two
sequences of the segments of billiard trajectories within the
ellipsoid $\mathcal E$ in $\mathbf E^d$, sharing the same $d-1$
caustics. Suppose the pair $(a_0, b_0)$ is $s$-skew, and that by the
sequence of reflections on quadrics $\mathcal Q^1,\dots,\mathcal
Q^{s+1}$ the minimal billiard trajectory connecting $a_0$ to $b_0$
is realized.

\smallskip

Then, each pair $(a_m, b_m)$ is $s$-skew, and the minimal billiard
trajectory connecting these two lines is determined by the sequence
of reflections on the same quadrics $\mathcal Q^1,\dots,\mathcal
Q^{s+1}$.
\end{theorem}

\begin{proof}
This may be proved by the use of the Double Reflection Theorem,
similarly as in Theorem \ref{th:gen.grida}.
\end{proof}

This theorem also can be stated for an arbitrary quadric, not only
for ellipsoids.

\begin{proposition}
Let $\mathcal E$ be an ellipsoid in $\mathbf E^d$,
 $\mathcal Q_1,\dots,\mathcal Q_{d-1}$
quadrics confocal to $\mathcal E$, and $k$ an integer.

\smallskip

Suppose that there exists a trajectory $(\ell_m)$ of the billiard
within $\mathcal E$, having the caustics
 $\mathcal Q_1,\dots,\mathcal Q_{d-1}$,
such that the pair $(\ell_0, \ell_k)$ is $s$-skew and that the
minimal billiard trajectory connecting $\ell_0$ to $\ell_k$ is
realized by the sequence of reflections on quadrics
 $\mathcal Q^1,\dots,\mathcal Q^{s+1}$
confocal to $\mathcal E$.

\smallskip

Then, for any billiard trajectory $(\hat\ell_m)$ within $\mathcal E$
with the caustics
 $\mathcal Q_1$, \dots, $\mathcal Q_{d-1}$,
the pairs of lines $(\hat\ell_m, \hat\ell_{m+k})$ are $s$-skew and
the minimal billiard trajectory between $\hat\ell_m$ and
$\hat\ell_{m+k}$ is realized by the sequence of reflections on
$\mathcal Q^1,\dots,\mathcal Q^{s+1}$.
\end{proposition}

\section{Conclusion}\label{sec:conclusion}

{\em ``One of the most important and also most beautiful theorems in
classical geometry is that of Poncelet (...) His proof was synthetic
and somewhat elaborate in what was to become the predominant style
in projective geometry of last century. Slightly thereafter, Jacobi
gave another argument based on the additional theorem for elliptic
functions. In fact, as will be seen below, the Poncelet theorem and
additional theorem are essentially equivalent, so that at least in
principle Poncelet gave a synthetic derivation of the group law on
an elliptic curve. Because of the appeal of the Poncelet theorem it
seems reasonable to look for higher-dimenisonal analogues...
Although this has not yet turned out to be the case in the
Poncelet-type problems...''}

\smallskip

These introductory words from \cite{GH1} written by Griffiths and
Harris exactly 30 years ago, could serve as a motto for the present
paper, announcing the programme realized here (see also
\cite{Pr,Pr2,Ves}).

\smallskip

Usually, the Poncelet theorem is associated with the points of
finite order on an elliptic curve. However, our analysis from the
Section \ref{sekcija:bilijarskaalgebra} shows that a Poncelet
polygon should be understood as a pair of equivalent divisors and
through a billiard analogue of the Abel theorem it can further be
associated to a meromorphic function of a hyperelliptic curve. On
the other hand, the billiard analogue of the Riemann-Roch theorem
gives uniqueness of the billiard trajectory with minimal number of
lines connecting two given lines in $\mathcal A_{\ell}$. Thus, the
Poncelet story in higher dimensions is at least as rich as a story
of meromorphic functions on hyperelliptic curves.

\smallskip

Such a deep relationship between hyperelliptic Jacobians and pencils
of quadrics through integrable billiard systems opens a new view to
the well-known, but still amazing role of elliptic coordinates in
the theory of integrable systems and in the separation of variables
in Hamilton-Jacobi method. The question of synthetic approach to the
addition law on nonhyperelliptic Jacobians remains open and could
lead to new methods of separation of variables.

\subsection*{Table of Notations}
\addtocontents{toc}{\numberline{}{\bf Table of Notations
\hfill\thepage}\newline}

\begin{itemize}

\item[$\mathcal A_{\ell}$\quad]
the family of all lines in $\mathbf E^d$ which are tangent to given
$d-1$ confocal quadrics (see Section \ref{sekcija:morfizam})

\item[$\mathcal C_{\ell},\mathcal C_{\mathcal O}$\quad]
the generalized Cayley curve (see Definition \ref{def:gen.Cayley} in
Section \ref{gen.lebeg})

\item[$\mathbf E^d$\quad]
the $d$-dimensional Euclidean space

\item[$F(V)$\quad]
the set of all $(d-2)$-dimensional linear subspaces of the
intersection $V$ of two quadrics in $\mathbf P^{2d-1}$ (see Lemma
\ref{lema:veza} in Section \ref{sekcija:morfizam})

\item[$\mathbf P^d$\quad]
the $d$-dimensional projective space

\item[$\mathbf P^{d*}$\quad]
the dual $d$-dimensional projective space

\item[$\mathcal Q$\quad]
a quadric

\item[$\mathcal Q^*$\quad]
the quadric dual to $\mathcal Q$

\item[$\mathcal Q_{\lambda}$\quad]
a quadric from the family (\ref{eq:confocal.family}) of confocal
quadrics in $\mathbf E^d$ (see Definition \ref{def:E.confocal} from
Section \ref{sec:confocal.E})
\end{itemize}

\subsection*{Acknowledgements}
\addtocontents{toc}{\numberline{}\newline}
\addtocontents{toc}{\numberline{}{\bf Acknowledgements
\hfill\thepage}\newline}

The research was partially supported by the Serbian Ministry of
Science and Technology, Project {\it Geometry and Topology of
Manifolds and Integrable Dynamical Systems}. One of the authors
(M.R.) acknowledges her gratitude to Prof.\ V.\ Rom-Kedar and the
Weizmann Institute of Science for the kind hospitality and support
during the work on this paper. The authors would like to thank the
referee for the valuable remarks and suggestions, which improved the
manuscript.

\end{document}